\documentclass[]{article}
\usepackage{amsfonts}
\usepackage{hyperref}
\usepackage[utf8]{inputenc}
\usepackage{listings}
\usepackage{tikz}
\usepackage{pgf}
\usepackage{siunitx}
\usetikzlibrary{positioning}
\usepackage[ruled, linesnumbered]{algorithm2e}
\usepackage{hhline}
\usepackage{tabu}
\usepackage{mathtools}
\DeclarePairedDelimiter{\ceil}{\lceil}{\rceil}
\DeclarePairedDelimiter{\floor}{\lfloor}{\rfloor}
\usepackage{svg}
\usepackage{comment}

\usepackage{color}
\definecolor{deepblue}{rgb}{0,0,0.5}
\definecolor{deepred}{rgb}{0,0,0.5}
\definecolor{deepgreen}{rgb}{0,0,0.5}

\newcommand\pythonstyle{\lstset{
	language=Python,
	basicstyle=\ttm,
	morekeywords={self},
	keywordstyle=\ttb\color{deepblue},
	emph={MyClass, __init__},
	emphstyle=\ttb\color{deepred},
	stringstyle=\color{deepgreen},
	frame=tb,
	showstringspaces=false}}

\lstnewenvironment{python}[1][]{
	\pythonstyle
	\lstset{#1}
	}
{}

\newcommand\pythoninline[1]{{\pythonstyle\lstinline!#1!}}
\title{Saving proof-of-work by hierarchical block structure\footnote{This is a republication of our work \cite{Somnium2022} under a pseudonym. For more details, see Section \ref{sec:republishing}}}
\author{Valdemar Melicher\\email: \href{mailto: valdemar.melicher@gmail.com}{valdemar.melicher@gmail.com}}

\newcommand{\pe}{e}
\newcommand{\hr}{h}
\newcommand{\thr}{h_B}
\newcommand{\ov}{v}
\newcommand{\vpb}{\beta}
\newcommand{\tec}{tec}
\newcommand{\cpt}{cpt}
\newcommand{\ep}{p}

\newcommand{\tf}{f}
\newcommand{\tfpb}{\phi}
\newcommand{\tr}{tr}
\newcommand{\epb}{E}
\newcommand{\wt}{w_B}
\newcommand{\state}{s}
\newcommand{\transf}{\tau}
\newcommand{\tx}{t}
\newcommand{\txerror}{e}
\newcommand{\etx}{\tx_e}

\newcommand{\fh}{\lambda}
\newcommand{\diff}{\delta}
\newcommand{\Diff}{\Delta}
\newcommand{\blocktime}{T}
\newcommand{\secc}{s}
\newcommand{\mc}{c}
\newcommand{\thrtocom}{\gamma}
\newcommand{\hasht}{t_h}
\newcommand{\ti}{\eta}
\newcommand{\bsize}{b}
\newcommand{\wtps}{\zeta}
\newcommand{\block}{\mathcal{B}}
\newcommand{\cti}{c_{\ti}}
\newcommand{\br}{r}
\newcommand{\shardf}{s}
\newcommand{\noc}{c}
\newcommand{\mfn}{{MFN}}

\newcommand{\NN}{\mathbb{N}}
\newcommand{\RR}{\mathbb{R}}
\newcommand{\set}[1]{\{#1\}}

\newcommand{\txset}{\mathcal{T}}

\newcommand{\notxs}{N}
\newcommand{\iset}{\mathcal{I}}
\newcommand{\isetl}{\mathcal{L}}

\newcommand{\Sum}[2]{\the\numexpr #1 + #2 \relax}

\newcommand{\refe}[1]{(\ref{#1})}

\def\bitcoinA{%
	\leavevmode
	\vtop{\offinterlineskip %\bfseries
		\setbox0=\hbox{B}%
		\setbox1=\hbox to\wd0{\hfil\hskip-0.03em
			\vrule height .3ex width .15ex\hskip .08em
			\vrule height .3ex width .15ex\hfil}
		\setbox2=\hbox to\wd0{\hfil\hskip-0.03em
			\vrule height .3ex width .15ex\hskip .08em
			\vrule height .3ex width .15ex\hfil}
		\vbox{\box1\box0}\box2}
}

\def\bitcoinAcorrected{%
	\leavevmode
	\vtop{\offinterlineskip %\bfseries
		\setbox0=\hbox{B}%
		\setbox1=\hbox to\wd0{\hfil\hskip-0.03em
		\vrule height .3ex width .15ex\hskip .08em
		\vrule height .3ex width .15ex\hfil}
		\setbox2=\hbox to\wd0{\hfil\hskip-0.03em
		\vrule height .3ex width .15ex\hskip .08em
		\vrule height .3ex width .15ex\hfil}
		\vbox{\box1\box0}\vskip-0.003em\box2}
	}
	
\def\bitcoinAd{%
	\leavevmode
	\vtop{\offinterlineskip %\textit
		\setbox0=\hbox{B}%
		\setbox1=\hbox to\wd0{\hfil\hskip0.27em
			\vrule height .3ex width .15ex\hskip .08em
			\vrule height .3ex width .15ex\hfil}%
		\setbox2=\hbox to\wd0{\hfil\hskip-0.03em
    		\vrule height .3ex width .15ex\hskip .08em
			\vrule height .3ex width .15ex\hfil}%
		\vbox{\box1\box0}\box2}
}	
	
\newtheorem{definition}{Definition}	
\newtheorem{hypothesis}{Hypothesis}	
\newtheorem{assumption}{Assumption}	

\newenvironment{algo}{}{}
\begin{document}

\maketitle

\begin{center}
\it To my wonderful wife, our children and yours as well
\end{center}

\begin{abstract}
	
We argue that the current \emph{Proof of Work} based consensus algorithm of the Bitcoin network suffers from a fundamental economic discrepancy between the real-world transaction costs incurred by miners and the wealth that is being transacted. Put simply, whether one transacts 1 satoshi or 1 bitcoin, the same amount of electricity is needed when including this transaction into a block. The notorious Bitcoin blockchain problems such as its high energy usage per transaction or its scalability issues are, either partially or fully, mere consequences of this fundamental economic inconsistency. We propose making the computational cost of securing the transactions proportional to the wealth being transfered, at least temporarily. 

First, we present a simple incentive based model of Bitcoin's security. Then, guided by this model, we augment each transaction by two parameters, one controlling the time spent securing this transaction and the second determining the fraction of the network used to accomplish this. The current Bitcoin transactions are naturally embedded into this parametrized space. Then we introduce a sequence of hierarchical block structures (HBSs) containing these parametrized transactions. The first of those HBSs exploits only a single degree of freedom of the extended transaction, namely the time investment, but it allows already for transactions with a variable level of trust together with aligned network fees and energy usage. In principle, the last HBS should scale to tens of thousands timely transactions per second while preserving what the previous HBSs achieved. 

We also propose a simple homotopy based transition mechanism which enables us to  relatively safely and continuously introduce new HBSs into the existing blockchain. 

Our approach is constructive and as rigorous as possible and we attempt to analyze all aspects of these developments, al least at a conceptual level. The process is supported by evaluation on recent transaction data.

\end{abstract}

\section{Introduction}
\label{sec:introduction}
Arguably, in 2008 a new era of finance began. In a groundbreaking paper \cite{Nakamoto2008}, Satoshi Nakamoto described a pure peer-to-peer payment system of Bitcoin. He presented an elegant and rather simple way how to prevent double spending in decentralized distributed networks without a trusted party using a \emph{Proof of Work} (PoW) concept from \cite{Back2002}.

The Bitcoin payment system, despite being extremely well thought out and elegant, could not grow up, at least yet, to the very expectations of its creators. 

Let us read Satoshi Nakamoto's statement from Section $1$ of \cite{Nakamoto2008}, arguing against financial institution mediation: \emph{``\dots The cost of mediation increases transaction costs, limiting the minimum practical transaction size and cutting off the possibility for small casual transactions,\dots''}. 

Unfortunately, looking at the current transaction costs on the main Bitcoin network or its available transaction throughput, it is obvious that Bitcoin suffers from the same symptoms but caused by a different illness. Let us paraphrase the Bitcoin scalability problem in a sentence similar to Satoshi Nakamoto's one: \emph{The costs of requiring that all transactions should be equal and all of them seen by all the nodes in the network increases transaction costs, limiting the transaction throughput and cutting off the possibility for small casual transactions.}

The limited transaction throughput is arguably the primary technological hurdle of Bitcoin's use as means of payment. And the high transaction costs are the main economical hurdle. But as indicated above, these are mere effects of uniformity and inflexibility of Bitcoin transactions, which do not respect economic realities and consequently cannot map to many real world situations. Put simply, if we require that each transaction is seen, signed and stored by everybody, the system has to be costly, slow and almost impossible to scale.

This paper is about, to our best knowledge, an innovative path of possible practical extensions of Bitcoin's PoW  based consensus algorithm and its protocol, starting form the definition of its transaction. We need to relax and expand the notion of a Bitcoin transaction and also of its block and network structure, all those aspects, in order to resolve its scalability issues. 

In Section \ref{sec:motivation} we first discuss the concepts of trust and risk, focusing on their multi-scale nature, which is in strong contrast to flat security guaranties and costs of the current Bitcoin implementation. This fundamental discrepancy essentially locks the Bitcoin development in a suboptimal state and prevents further innovation. In the followings sections we gradually annihilate this inconsistency, while preserving Bitcoin's simplicity to a maximal possible degree. Importantly, transactions would still be saved in an adapted multi-layered blockchain and secured by the mainnet and not a distinctively separate second layer technology such as \emph{Lightning Network} \cite{Beres2020}.

In preparatory Section \ref{sec:par_trans}, we discuss security, fees and electric power consumption of Bitcoin transactions and the alignment of these three aspects.  We conclude the section by definition of a 2-parameter family of extended transactions.
In Section \ref{sec:hbs} we propose a sequence of hierarchical block structures (HBSs) containing those parametrized transactions. Starting from the first HBS presented in Section \ref{sec:hbs_svc}, users can choose a level of security which they require to transact and this in a transparent and simple way to grasp. And the costs will be proportional to this level of trust. Sections \ref{sec:hbs_binary_tree} and \ref{sec:hbs_concurrent_levels} are devoted to the scaling of on-chain transaction capacity. Some supporting material is presented in Appendix \ref{app:bitcoin_economics} and cited when needed. Particularly, alternative consensus algorithms are shortly discussed in Section \ref{sec:current_state_blockchain_affairs}. For reproducibility of results, almost all algorithms used in this publication are published on github \cite{github} and for completeness also presented in Appendix \ref{app:algorithms}. 

Finally, let us note that while the ideas that are presented in this paper are certainly widely applicable outside of Bitcoin's realm, their possible applications in different contexts are left for the reader.

\section{Rigidity of Bitcoin's trust model}
\label{sec:motivation}

\subsection{Distributed trust model avoids systemic risks}
Bitcoin allows for direct peer-to-peer online payments without an all-knowing trusted party(ies). In essence, the trust is distributed in the network among its participants, particularly among its mining nodes. The information exchange among the participants is governed by Bitcoin's protocol and its rules. Due to its distributed nature, Bitcoin avoids all kinds of systemic risks, e.g. it cannot be easily globally censored and it does not have a single point of failure and it is permissionless, it is anyone can join and leave network at any time.

\subsection{At individual level risks are countless}
All sorts of other risks remain. For example, the Bitcoin offers irreversible transactions. This is an advantage for sellers, not so much for buyers. Satoshi Nakamoto, being aware of it, proposed routine escrow mechanisms to protect buyers \cite{Nakamoto2008}. But an escrow agent (e.g. implemented via multisig) is a trusted third party. Or imagine that you buy a property and pay in bitcoins. The ownership of the property and/or the coins can be disputed by either party or even a third party at any moment, i.e. the parties are exposed to legal risk and the judiciary might become a trusted party. Further, we all are susceptible to $5\$$ wrench attacks. Last but not least, humans are fallible and we often loose access to private keys. 

The above mentioned risks are mostly local, either transaction specific risks having to do with the fact that their underlying contracts relate to values in physical world or any other risks at individual level having to do with our sheer existence in real world and our own fallibility.

Altogether, as the scale goes from global to local, both trust and risks become concentrated when making transactions. As individuals, we have to be and are prepared to face a myriad of foreseen and unforeseen circumstances. Based on risks we face and our abilities to mitigate them, trust we have in our counterparty and taking into account transaction costs and rewards, we either decide to transact or not. 

\subsection{Relaxing trust model of consensus algorithm?}
\label{sec:relax_model_of_consensus}
The main recurrent costs for bitcoin miners are electric energy bills, see Appendix \ref{app:mining_economics}. Since the hashing time is linearly proportional to size of data being hashed, see Algorithm \ref{alg:hashing_time}, it seems natural to compute fees based on transaction size, as it is currently done. And indeed, purely from the miners' perspective, whether a transaction moves 100 BTC or 1 Satoshi, it does not matter if their bit size is the same. Is this not a problem? Should not moving a fortune cost more then paying few cents?

In our opinion, it clearly should. Otherwise, small payments will become clearly impossible as Bitcoin transitions from issuance system where miners are paid mainly by rewards to a system based on transactions fees, and this regardless whether they could be included into a block given the size constraints or not. The second but not secondary issue is that currently those transacting smaller amounts actually pay fees for security of those transacting bigger ones, promoting a rising inequality.

We think that this discrepancy between transaction fees and wealth being transfered is probably the most fundamental problem, stopping the Bitcoin network from becoming a general payment system as envisioned by Satoshi Nakamoto. So, how to cut the Gordian knot?

A solution becomes almost evident, if one aims for a \emph{transaction fee per byte} $\tfpb$ proportional to the wealth being transfered $\ov$. Then, at least conceptually, arbitrarily small payments become possible. And it clearly makes system at once more honest. A direct consequence of this aim is that the security of transactions has to be proportional to the wealth being transfered as well, at least temporarily. Otherwise a difficulty would arise how to pay the miners for their work.

And indeed, from the user risk/reward perspective described in the previous section is Bitcoin's current consensus protocol rather too coarse, idealistic and rigid. It simply does not allow a user to choose an appropriate costs/risk ratio. The transactions are public and every full node has to be aware of all the transactions ever included in the blockchain. Either a user accepts this absolute security and pays all the mining nodes for their work the corresponding price and bears the consequences (e.g. a long confirmation time) or he opts to not transact. The second option is a rational choice in many practical circumstances. Again, a successful general online payment system has to allow for inexpensive quick micro payments.

In the next sections, we will include \emph{level of trust in the payment system} as one of the risks and \emph{network fee} as one of the costs in individual's decision space, with some safe defaults which are keeping wealth being transfered and safety proportional.

\section{Extending transaction space}
\label{sec:par_trans}

The Bitcoin system can be seen as a \emph{state-transition machine} \cite{Buterin2013}. Its \emph{state} $\state$ consists essentially of all ownership records of all existing coins and it can be changed only according to the Bitcoin protocol by applying the corresponding \emph{state transition function} $\transf_B$. This function is a very handy abstraction, which allows one to effectively hide whole the Bitcoin protocol complexity under one symbol. It takes old state $\state_o$, a transaction $\tx$ and evaluates as
\begin{equation}
\transf_B(\state_o, \tx) = \left\{\begin{array}{l}
\state_n \quad\mbox{if $\tx$ is compliant},\\
(\state_o, \txerror ) \quad{\mbox{if $\tx$ is not compliant,}} 
\end{array}\right.
\label{eq:transition_function}
\end{equation}
where $\state_n$ is the new state and $\txerror$ is a returned error in the case of noncompliant $\tx$. Any transaction $\tx$ is essentially a set of instructions created by the owner of unspent bitcoins, called \emph{Unspent Transaction Outputs} (UTXOs), how to allocate those coins. As the name clearly suggests, UTXOs are outputs of previous transactions that have not yet been inputs of any already existing valid transaction. Each full node computes its UTXO set independently by validating all transactions in the blockchain. It recurrently applies Equation \refe{eq:transition_function} starting from bitcoins created as mining rewards in a \emph{coinbase} transaction.  Consequently, UTXOs are a derived representation of a current state $\state$ which is actually recorded in the blocks.

The main goal of this section is to define a concept of an extended Bitcoin transaction. Before we actually do so in Section \ref{sec:ext_tx}, we describe a simple incentive based economic model of Bitcoin transactions, focusing first on the aspect of security in Section \ref{sec:transaction_security} and then discussing transaction fees in Section \ref{sec:fees}.

\subsection{Transaction security}
\label{sec:transaction_security}
As explained in Appendix \ref{app:centralization_of_mining}, the long term investments of miners to acquire highly specialized non-repurposable hashing hardware can be seen as their main incentive to act honestly and keep the Bitcoin network running smoothly \cite{Hasu2019}. Any loss of trust would negatively impact the BTC market prize, decreasing the future returns of the miners.

Consequently, the security of Bitcoin transactions is proportional to \emph{miner commitment} $\mc[\bitcoinA]$, which in turn is proportional to the total hashrate of the network $\thr$ in $H/s$. Under the assumption of linearity we can stipulate that:\footnote{Occam's razor.}
\begin{equation}
c[\bitcoinA] : = \thrtocom [\mbox{\bitcoinA s/H}]\cdot \thr[\mbox{H/s}],
\label{eq:thrtocom}
\end{equation}
where $\gamma > 0$ is for simplicity assumed being a real constant. Since we do not know its value, the equation \refe{eq:thrtocom} is only qualitatively interesting. E.g. the dimension $\bitcoinA s/H$ of $\gamma$ tells us that the miners are committed, if their BTC gains per one hash are high and the ``time horizon'' of value extraction is long. From now on, the square brackets in formulas are exclusively used to present dimensions of variables, with an exception of a closed interval.

The miners have to find a nonce for their block such that its hash is smaller than a predefined target determined by the \emph{difficulty} $\Diff$ \cite{difficulty}. The difficulty is adjusted every 2016 blocks, i.e. approximately every 2 weeks, with the goal to achieve average block time of 10 minutes. Time, the difficulty $\Diff$ and the total Bitcoin network hashrate $\thr$ fulfill this equation
\begin{equation}
\thr\left[\frac{\mbox{H}}{\mbox{s}}\right] = \frac{1}{65535}\frac{2^{48}\Diff[\mbox{H}]}{600\mbox{s}} \approx  \frac{2^{32}\Diff[\mbox{H}]}{600\mbox{s}}.
\label{eq:difficulty}
\end{equation} 

The difficulty $\Diff$ is usually presented as being dimensionless quantity, but it actually represents the number of hashes computed by all Bitcoin miners in one block time, on average, divided by $2^{32}.$

Just another ingredient of the ``security soup'' is the realization that \emph{hashing time} $\hasht$ is linearly proportional to bitsize $\bsize(m)$ of a hashed message $m$, see Algorithm \ref{alg:hashing_time}. We can write:
\begin{equation}
\label{eq:hash_time}
\hasht(m)[\mbox{s}] := \ti\left[\frac{\mbox{s}}{\mbox{b}}\right] \cdot \bsize(m)[\mbox{b}], 
\end{equation}
where variable $\ti$ represents the time necessary to sign one bit of information. Let us call it \emph{time investment}. In the case of signing by all Bitcoin miners, the corresponding $\ti_B$ can be estimated easily from block size $\bsize(\block)$ in bits:
\begin{equation}
\label{eq:eta_B_estimate}
\hat\ti_B\left[\frac{\mbox{s}}{\mbox{b}}\right] = \frac{600\mbox{s}}{\overline{\bsize(\block)}[\mbox{b}]}, 
\end{equation}
where $\overline{\bsize(\block)}$ depicts the average block size and subscript $B$ means we consider the whole Bitcoin network. From now on, $\overline{x}$ will depict the sample mean value of $x$ and $\hat{x}$ a general estimator of $x$.
From block size limit of e.g. $1$MB, we obtain a lower bound on $\ti_B$ of approximately $\num{7.15e-5}.$

We propose the following simple incentive based Bitcoin security model: 
\begin{definition}
	\label{def:security}
	A security $\secc(t)$ of transaction $t$ with output value $\ov(t)$ is defined as	
	\begin{equation}
	\secc(t)\left[\frac{H}{\bitcoinAd}\right] := \frac{\hr(t)[H/s]\cdot\hasht(t)[s]}{\ov(t)[\bitcoinAd]},
	\label{eq:def_security}
	\end{equation}
	where $\hr(t)[H/s]$ is the hashrate used to sign $t$ in time $\hasht(t)$.
\end{definition}
This definition is intuitively very easy to understand. The security is proportional to the hashing resources used and to the time they are employed to compute the transaction signature. And it is inversely proportional to the value transfered since it is the eventual gain from the double spending attack targeting this transaction.

Using \refe{eq:hash_time} we obtain for the security $\secc_B$ of current Bitcoin transactions:
\begin{equation}
\secc_B(t) = \frac{\thr\cdot\ti_B\cdot\bsize(t) }{\ov(t)}.
\label{eq:def_security_B}
\end{equation}
which leads to an average per block $\block$ security
\begin{equation}
\overline{\secc_B(t)}|_{t\in\block} = \frac{\thr\cdot\ti_B}{|\block|}\sum_{t\in \block}\frac{\bsize(t)}{\ov(t)},
\label{eq:def_security_average}
\end{equation}
where $|\block|$ represents the cardinality of set $\block$, i.e. here the number of transactions in $\block$.

Definition \ref{def:security} describes the security from the point of view of individual transactions. Currently, transactions with diverse output values $\ov$ are included into the same block. Following the developments, security per block $\block$ can be consistently defined as:
\begin{equation}
\secc_B(\block) := \frac{\thr\cdot \hasht(\block)}{\ov(\block)} = \frac{\thr\cdot\ti_B\cdot\bsize(\block) }{\ov(\block)}
\mathop{\approx}\limits^{\refe{eq:eta_B_estimate}}
\frac{\thr\cdot 600s}{\ov(\block)}\cdot\frac{\bsize(\block)}{\overline{\bsize(\block)}}
\mathop{\approx}\limits^{\refe{eq:difficulty}} \frac{2^{32}\Diff}{\ov(\block)}
\label{eq:def_security_per_block},
\end{equation}
which sheds more light on relation between security and difficulty. The numerator of \refe{eq:def_security_per_block} is an equivalent of difficulty, up to multiplication by a constant. When the total output value $\ov(\block)$ per block $\block$ rises in time, difficulty $\Diff$ has to rise at least proportionally to preserve security. Based on the alignment of \refe{eq:def_security_per_block} with \refe{eq:def_security}, we can define the transaction difficulty $\diff(t)$ of a transaction $t$ accordingly:
\begin{equation}
\diff(t)[\mbox{H}] := \hr(t)\cdot\hasht(t) \mathop{=}\limits^{\refe{eq:hash_time}} \hr(t)\cdot\ti(t)\cdot\bsize(t)
\label{def:tx_diff}
\end{equation}
The current security per block $\secc_B(\block)$ can be estimated as
\begin{equation}
\hat\secc_B(\block) = \frac{2^{32}\Diff}{\overline{\ov(\block)}}
\label{eq:def_security_per_block_estimate},
\end{equation}
where average output value per block $\overline{\ov(\block)}$ can be computed from real blockchain data, see e.g. Section \ref{sec:segmentation_of_transactions}.

Further, we can write \refe{eq:def_security_per_block} as
\begin{equation}
\secc_B(\block) =  \frac{\thr\cdot\ti_B\cdot\sum_{t\in\block}\bsize(t) }{\sum_{t\in\block}\ov(t)}=\frac{\thr\cdot\ti_B\cdot\overline{\bsize(t)}}{\overline{\ov(t)}},
\label{eq:def_security_per_block_average_tx}
\end{equation}
i.e. in terms of the average transaction where the sample is the whole block $\block$. For a message with average size and output values, we can conclude by comparing \refe{eq:def_security} and \refe{eq:def_security_per_block_average_tx} that the current block-based bitcoin security is identical to the individual one according to Definition \ref{def:security}.

However, if e.g. an average sized transaction has a much higher output value than the average one is, the current block-based security may be insufficient. On the contrary, if e.g. an average sized transaction has a much lower output value than the average one is, this security may be an overkill. Let us make the following assumption:
\begin{assumption}
	\label{ass:same_security}
	If one transacts, he expects the same security $\secc(t)$ regardless of transaction value $\ov(t)$.
\end{assumption} 
This formally under the validity of model 	\refe{eq:def_security}  means that there exists a $c\in\RR^+$ such that transaction difficulty $\diff(t) = c\cdot\ov(t)$ for all admissible transactions $t,$ implying $\secc(t) = c.$ This $c$ is certainly at least time dependent, but it is independent on individual transactions $t$.

\subsection{Transaction fees}
\label{sec:fees}
Currently, Bitcoin users pay fees that are only proportional to transaction size $\bsize(t)$. From the perspective of miners, this fee scheme is an optimal one. Their main running costs are payments for electric power consumption used to hash blocks of transactions. In Equation \refe{eq:tec} all the other variables except $\bsize(t)$ are transaction invariant. First, all the transactions are included in a block hashed by all the miners with their combined hash rate $\hr_B$. Second, blocks are created on average every 10 minutes. 

However from the usefulness perspective, as already informally argued in Section \ref{sec:relax_model_of_consensus}, the current fee scheme based only on $\bsize(t)$ is far from optimal. We will look for an optimal fee scheme(s) from both user and miner perspectives in a relaxed setting where strict Bitcoin constraints on hash rate and hashing time are removed.

In \refe{def:tx_diff}, the definition of transaction difficulty $\diff(t)$, both the hashing time $\hasht(t)=\ti(t)\cdot\bsize(t)$ and hash rate $\hr(t)$ are transaction dependent variables by design.  Following \refe{eq:tec}, $\diff(t)$ multiplied by energy efficiency $\pe$ of mining hardware gives us energy consumption of a transaction. Thus in any case, to fairly cover miner costs, users should  pay transaction fees $\tf$ proportional to $\diff(t)$:
\begin{equation}
\tf [\bitcoinA] \propto \diff(t)[\mbox{H}].
\label{eq:fee_prop_to_diff}
\end{equation}
In effect,  users then still pay fees proportional to transaction size $\bsize(t)$ but also proportional to variable hashing capacity used $\hr(t)$ and its employment time, determined by $\ti(t).$

The transaction difficulty $\diff(t)$ is the numerator of \refe{eq:def_security}, which defines the security of individual transactions $\secc(t).$  Thus $\secc(t)$, except being proportional to the long term miner commitment $\mc$ through $\hr(t)$, is also proportional to the energy used to sign the transaction, i.e. to the main running cost of miners, making it a rather robust concept. The more electrical energy is used to sign a transaction, the more of it is needed to double spent it.

%In the end, the Bitcoin system is a public ledger of transactions that are signed by the hashing power of miners and that is what the users should be paying for, both for its capacity and employment time based on \refe{def:tx_diff}
Finally, if we want to adapt the Bitcoin ecosystem to be conformant to Assumption \ref{ass:same_security}, we have to find mappings $\hasht(t)=\ti(t)\cdot\bsize(t)$ and $\hr(t)$ and constant $c\in\RR^+$ such that transaction difficulty $\diff(t) = \hasht(t)\cdot\hr(t) = c\cdot\ov(t)$ for all admissible transactions. Assumption \ref{ass:same_security} thus implies that the transaction fees $\tf$ should be always proportional to the transaction value $\ov(t):$
\begin{equation}
\tf [\bitcoinA] \propto \ov(t)[\bitcoinA].
\end{equation}

\subsection{Definition of the extended transaction}
\label{sec:ext_tx}

Now, we are ready to introduce our main vehicle of change. Let us extend the Bitcoin transaction $\tx$ by two-parameters:

\begin{definition}
%$\etx(\dots, \etxd, \etxf)$ is an extended Bitcoin transaction, where parameter $\etxd\in(0,\diff]$ determines the relative difficulty mining blocks to which this transaction belongs, $\diff$ being the current difficulty of the whole Bitcoin network \footnote{Leading to an average block time of approximately $10$ minutes.} and $\etxf\in(0,1]$ determines the fraction of mining nodes which include this transaction into their blocks.
$\etx(\dots, \fh, \ti)$ is an extended Bitcoin transaction, where $\fh\in(0,1]$ determines the fraction of all mining nodes which include this transaction into their blocks and the time investment $\ti$ defined in \refe{eq:hash_time} determines how long these miners on average look for an admissible hash of this transaction.
\label{def:extended_transaction}
\end{definition}

The current Bictoin transactions $t$ are naturally embedded in the extended transaction ``space''. We have $\etx(\dots, 1, \ti_B) = t(\dots)$, i.e. the current Bitcoin transactions are included into blocks by all the available mining nodes and it takes the nodes on average 10 minutes to sign the blocks with the time investment $\ti_B$ estimated in \refe{eq:eta_B_estimate}. This can be immediately expressed in difficulty language using \refe{eq:difficulty}:
\begin{equation}
\ti_B \approx \frac{2^{32}\Diff}{\hr_B \cdot \overline{\bsize(\block)}}
\end{equation}
How to understand Definition \ref{def:extended_transaction} in a situation when $\tx \neq \etx$? For example $\etx(\dots, 0.5, \ti_B)$ is a transaction mined by only a half of all mining nodes and the hashing difficulty is reduced by factor two which essentially means that the resulting hash can have a one less leading zero since
\[
\ti_B \approx \frac{2^{32}(\Diff/2)}{(\hr_B/2) \cdot \overline{\bsize(\block)}}.
\]
On average, a block containing such transactions will be created in the same time as current Bitcoin blocks under the constraint that the block size does not change.

Two extra parameters of the extended transaction $\etx$ control two degrees of freedom of its difficulty defined in \refe{def:tx_diff}:
\begin{equation}
\diff(\etx) = \hr(\etx)\cdot\ti(\etx)\cdot\bsize(\etx) = \fh(\etx)\hr_B\cdot\ti(\etx)\cdot\bsize(\etx),
\label{eq:difficulty_control}
\end{equation}
i.e. the hashing resources used and the time invested to sign $\etx$.

We would like to recall that both the security $\secc(\etx)$ and the transaction fees $\tf(\etx)$ are quite naturally proportional to $\diff(\etx),$ see \refe{eq:def_security} and \refe{eq:fee_prop_to_diff} respectively. Thus, by controlling $\fh$ and $\ti$ we control transaction security and at the same time offer transaction fees proportional to this security. This is indeed not by coincidence but by respecting the Bitcoin economic realities while deriving the Bitcoin incentive based model presented in the previous sections. It also means that the model does not suffer from some obvious inconsistencies. From now on we use notation $t$ instead of $\etx$ since we will only consider extended transactions.

\section{Hierarchical block structures (HBSs)}
\label{sec:hbs}
In the subsequent sections,  we will employ the concept of extended transaction $t$ from Definition \ref{def:extended_transaction} to suggest modifications of the Bitcoin block structure that allow users to transact with different levels of security, while aiming costs proportional to this security. Of course, users will be free to override any defaults and to transact using an increased level of security if they are prepared to pay higher fees or to decrease the default security, if they feel that a lower one is sufficient and thus pay even less. Based on freedom of choice, a ``marketplace of trust'' should emerge quite naturally.      

The first proposition in Section \ref{sec:hbs_svc} would require minimal changes to the current Bitcoin ecosystem. Its further evolutions in Sections \ref{sec:hbs_binary_tree} and \ref{sec:hbs_concurrent_levels} lead to both more usability, complexity and more disruption of Bitcoin's status quo.

\subsection{Single virtual computer}
\label{sec:hbs_svc}

In this section we still assume that all transactions are seen by all mining and full nodes, i.e. Bitcoin network remains acting as a single ``homogeneous'' state-transition machine. Consequently, the scalability issues will remain unresolved. But we achieve the first necessary condition for Bitcoin to potentially become a successful general on-chain payment system: the ability to transact at different security levels and pay proportional fees.

We consider only transactions $t$ signed by all miners, i.e. $\fh(t) = 1$ for any admissible extended transaction $t,$ which trivially implies $\hr(t)=\fh(t)\hr_B = \hr_B.$ The only degree of freedom of the extended transaction $t$ we will tune here is the time investment $\ti.$

Let $\blocktime$ denote the current average block time resulting from the difficulty $\Diff$. In this section, instead of mining just one block in $\blocktime$, we propose mining a sequence of blocks containing transactions with different levels of security. As depicted in Figure \ref{fig:hbs_svc_block_scheme}, the blocks with $L\in\NN^+$ different security levels governed by time investments $\ti_l, l = 0, \dots, L-1$ are periodically repeating, forming imaginary super-blocks. In each super-block the lowest block contains transactions which require a minimal security $\ti_{L-1}$ and the highest block the transactions which require the maximal security $\ti_0$ and $\ti$ is monotonically decreasing in $l$, i.e. 
\begin{equation}
\ti_{l+1} < \ti_l \quad\forall l \in \{0,\dots, L-2\}.
\label{eq:monotonicity_of_ti}
\end{equation}

\newcommand{\hsbFigOneShift}{4}
\begin{figure}
	\hspace{-1cm}
	\begin{tikzpicture}
	\node[draw=none] at (0,4)	(level1) {$\phantom{\qquad\qquad\quad}\dots$};

	\node[draw] at (\hsbFigOneShift,4)	(level1_a) {$\{t(\dots, 1, \ti_0)\}$};
    \node[draw] at (\hsbFigOneShift,3)	(level2_a) {$\{t(\dots, 1, \ti_1)\}$};
    \node[draw] at (\hsbFigOneShift,2)	(level3_a) {$\{t(\dots, 1, \ti_2)\}$};
	\node[] at (\hsbFigOneShift,1)    (levelx_a) {$\dots$};
	\node[draw] at (\hsbFigOneShift,0)	(levell_a) {$\{t(\dots, 1, \ti_{L-1})\}$};
	\draw[-stealth] (levell_a) -- (levelx_a);
	\draw[-stealth] (levelx_a) -- (level3_a);
	\draw[-stealth] (level3_a) -- (level2_a);
	\draw[-stealth] (level2_a) -- (level1_a);
	
	\draw[-stealth] (level1.east) -- (levell_a.west);
	
	\node[draw] at (2*\hsbFigOneShift,4)	(level1_b) {$\{t(\dots, 1, \ti_0)\}$};
	\node[draw] at (2*\hsbFigOneShift,3)	(level2_b) {$\{t(\dots, 1, \ti_1)\}$};
	\node[draw] at (2*\hsbFigOneShift,2)	(level3_b) {$\{t(\dots, 1, \ti_2)\}$};
	\node[] at (2*\hsbFigOneShift,1)    (levelx_b) {$\dots$};
	\node[draw] at (2*\hsbFigOneShift,0)	(levell_b) {$\{t(\dots, 1, \ti_{L-1})\}$};
	\draw[-stealth] (levell_b) -- (levelx_b);
	\draw[-stealth] (levelx_b) -- (level3_b);
	\draw[-stealth] (level3_b) -- (level2_b);
	\draw[-stealth] (level2_b) -- (level1_b);
	\draw[-stealth] (level1_a.east) -- (levell_b.west);

	\node[draw=none] at (3*\hsbFigOneShift,0)	(levell_c) {$\dots$\phantom{\qquad\qquad\quad}};

	\draw[-stealth] (level1_b.east) -- (levell_c.west);
	
	\draw (2.2,-0.5) -- (2.2,4.5);
	\draw (2.2,4.5) -- (5.8,4.5);
	\draw (5.8,4.5) -- (5.8, -0.5);
	\draw (5.8,-0.5) -- (2.2, -0.5);
	
	\draw (\hsbFigOneShift+2.2,-0.5) -- (\hsbFigOneShift+2.2,4.5);
	\draw (\hsbFigOneShift+2.2,4.5) -- (\hsbFigOneShift+5.8,4.5);
	\draw (\hsbFigOneShift+5.8,4.5) -- (\hsbFigOneShift+5.8, -0.5);
	\draw (\hsbFigOneShift+5.8,-0.5) -- (\hsbFigOneShift+2.2, -0.5);
	\end{tikzpicture}
	\caption{Blockchain of blocks with different security levels. Each sub-block represents a current Bitcoin block with a different hashing difficulty.}
	\label{fig:hbs_svc_block_scheme}
\end{figure}
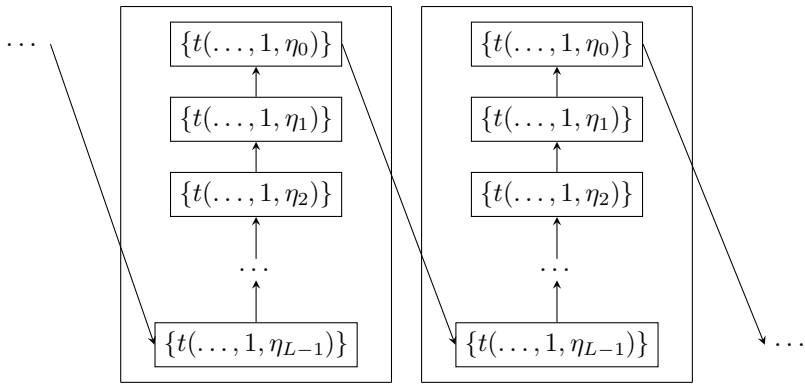

\subsubsection{Desing methodology}
We approach the design of our first Bitcoin evolution prototype in two steps:
\begin{enumerate}
\item We segment Bitcoin transactions waiting in the mempool $\txset$ for confirmation according to their output \emph{values per bit} $\vpb$, i.e. we handle the ratio
\begin{equation}
\vpb(t) :=\frac{\ov(t)}{\bsize(t)}
\end{equation}
occurring in  \refe{eq:def_security}. Then we will assign them into individual security level sub-blocks from Figure \ref{fig:hbs_svc_block_scheme}.
\item Then we determine variable $\ti_l$ for each level $l = 0, \dots, L-1$ such that security complies with Assumption \ref{ass:same_security}.
\end{enumerate}

There certainly exists no single best way how to accomplish those tasks but we will present our best candidate and it is up to the Bitcoin community to eventually refine it, fill in details or alternatively propose a different solution altogether.

The developments are guided by performance on recent transaction data of the Bitcoin network, also because any eventual candidate for next Bitcoin protocol has to allow for smooth transition from the current state. One can employ Algorithm \ref{alg:tx_dataset} from Appendix \ref{app:algorithms} to create a random sample of approximately 500 recent blocks starting from block height $650000$ and ending with $700000$. Our set is published on github \cite{github}, where the reader also finds the majority of code used for this article.

\subsubsection{Segmentation of transactions}
\label{sec:segmentation_of_transactions}

\begin{figure}[ht]
	\scalebox{0.85}{\input{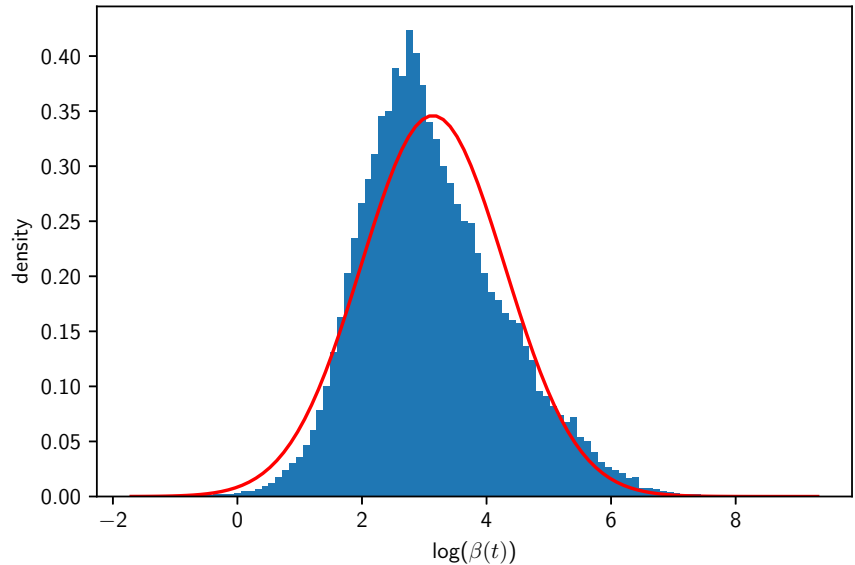}}
	\caption{Distribution of per transaction aggregated output values plus the corresponding log-normal distribution}
	\label{fig:tx_vpb}
\end{figure}

To illustrate the current distribution of aggregated output transaction values per bit $\vpb(t)$ of Bitcoin network, we employed Algorithm \ref{alg:plot_tx_hist} to plot a histogram of 
\[\lg\left(\vpb(t)\right),\]
where $\lg(\cdot) := \log_{10}(\cdot)$, together with the fitted log-normal probability density function, this resulting in Figure \ref{fig:tx_vpb}.

The distribution of $\vpb(t)$ is certainly not perfectly log-normal, but rather skewed to the right side with a long tail. One can only speculate why but two factors seem to be in play. First, as already explained, Bitcoin is not usable for everyday transactions, certainly not on chain. Second, high $\beta(t) $ transactions do not pay fair fees to miners for their work \footnote{The author is not involved in mining business.}. I.e. low-value transactions are discouraged (or impossible) and high-value ones encouraged. E.g. if a whale moves $10000\bitcoinA$ to a cold storage, it can cost less than $10^{-6}\%$ to do so. Nevertheless, the log-normal parametrization is already sufficiently accurate and the distribution of $\vpb(t)$ should become more log-normal and symmetric if Bitcoin becomes more useful for everyday payments on chain.

Altogether, we propose segmenting the unconfirmed Bitcoin transactions waiting in the mempool $\txset$ into $L\in\NN^+$ lists corresponding to individual levels from Figure \ref{fig:hbs_svc_block_scheme} based on their $\lg\left(\vpb(t)\right)$ values. The following simple algorithm can be employed:

\begin{algorithm}[H]
	\caption{A candidate for segmentation of a transaction set into subsets.}
	\label{alg:tx_segmentation}
	\SetAlgoLined
	\KwData{$L, N \in \NN^+$; a set of transaction tuples $\txset = \set{(t_i, \vpb_i)}_{i\in\iset}$, $\iset = \set{i| i = 0, \dots, \notxs-1}$, where $t_i$ are some unique transaction identifiers and $\vpb_i$ the corresponding transacted values per bit.}
	\KwResult{$L$ sets $\txset_i$, such that $\txset = \cup \txset_i$ and $\forall i,j\in\iset: \txset_i\cap\txset_j = \emptyset $}
	$\txset$ = sort set $\txset$ in descending order of $\vpb_i$\;
   	\For{$i\gets 0$ \KwTo $L-1$}{
   		$\txset_i = \{\}$\;
   	}
    $S_L = (\lg(\vpb_0)-\lg(\vpb_{N-1}))/L$ (or $S_L = (\ceil{\lg(\vpb_0))}-\floor{\lg(\vpb_{N-1}})/L$)\;
    $C = \lg(\vpb_0) - S_L$; $l = 0$\;
	\For{$i\gets 0$ \KwTo $N-1$}{
		\If{$\lg(\vpb_i) < C$}{
			$l = l+1$\;
		}
		$\txset_l = \txset_l \cup \{(t_i, \vpb_i)\}$\; 
	}
\end{algorithm}
The code is practically self-explanatory. The interval between the minimal $\lg(\vpb_{N-1})$ and the maximal $\lg(\vpb_0)$ is uniformly divided into $L$ subintervals and a transaction $t$ is assigned to set $\txset_0$ if
\[\lg(\vpb(t))\in [\lg(\vpb_0)-S_L, \lg(\vpb_0)]\]
and to $\txset_l$ for $l\in\left\{1,\dots, L-1\right\}$ if
\[\lg(\vpb(t))\in [\lg(\vpb_0)-(l+1)S_L, \lg(\vpb_0)-lS_L).\]

The transactions are thus segmented based only on $\vpb$ values and not on the frequency of their occurrence. We could instead compute L-quantiles of $\lg(\vpb(t))$ for our set of transactions $t$. Then all sets $\txset_l$ for $l\in\left\{0,\dots, L-1\right\}$ would have similar cardinalities. However, it is much easier for a user to understand mapping between the security levels and the transaction values/costs when the interval $[\lg(\vpb_0), \lg(\vpb_{L-1})]$ is uniformly divided. Also, this choice is more forward looking, which becomes clear in the next sections. 
%Nevertheless, it is just an essence... 

A Python implementation of Algorithm \ref{alg:tx_segmentation} is presented in Algorithm \ref{alg:tx_segmentation_python} in Appendix \ref{app:algorithms}. We employ it to segment our Bitcoin transaction dataset. The results, obtained by running Algorithm \ref{alg:get_table_segmented_txs}, are presented in Table \ref{table:segmented_txs}.

\def\arraystretch{1.1}
\begin{table}
	\begin{tabular}{|c|c|c|c|c|c|c|}
		\hline
		$l$ & 0 & 1 & 2 & 3 & 4 & 5\\ \hhline{|=|=|=|=|=|=|=|}
		$|\mathcal{T}_l|$ & 615 & 24967 & 220097 & 571652 & 130420 & 1861 \\ \hline
		min($\beta(t)$) & 3e+07 & 4.4e+05 & 6.4e+03 & 9.3e+01 & 1.4 & 0.02 \\ \hline
		max($\beta(t)$) & 2.1e+09 & 3e+07 & 4.4e+05 & 6.4e+03 & 9.3e+01 & 1.4 \\ \hline
		$\overline{\beta(t)}$ & 3.6e+08 & 2.3e+06 & 6.3e+04 & 1.3e+03 & 4.5e+01 & 0.73 \\ \hline
		min($v(t)$) & 5.4e+10 & 6.8e+08 & 9.6e+06 & 7.6e+04 & 2.1e+03 & 3.1e+02 \\ \hline
		max($v(t)$) & 6e+12 & 6.6e+11 & 1.6e+11 & 7.4e+09 & 1.2e+08 & 3e+05 \\ \hline
		$\overline{v(t)}$ & 8.7e+11 & 6.4e+09 & 2.3e+08 & 6.9e+06 & 2.6e+05 & 3.7e+03 \\ \hline
		$\sum v(t)$ & 5.3e+14 & 1.6e+14 & 5.1e+13 & 3.9e+12 & 3.4e+10 & 6.9e+06 \\ \hline
	\end{tabular}
	\caption{Distribution of transaction values after segmentation: $L=6$}
	\label{table:segmented_txs}
\end{table}

Let us discuss the results a little bit. We have divided our dataset containing almost $950k$ transactions from 493 blocks into only $6$ levels, so that Table \ref{table:segmented_txs} fits on page. One immediately recognizes how extreme this segmentation really is. 

Only $615$ transactions, it is less than $0.1\%$,  move more than $71\%$ of all the transacted wealth, which on average represents barely more than 1 transaction per block. However, together they only pay fees in the hundreds of USD.

On the other side of the spectrum in level $5$ are few, almost $2k$ transactions, which move on average 3700 sats, at current exchange rate $\mbox{BTCUSD}$ less then 2 USD. This is significantly less them our estimate of electricity costs per transaction of $14.3$ USD from Appendix \ref{app:energy}. We can argue based on this observation, that costs are a limiting factor for inclusion of transactions into the blockchain. And it is indeed a rational choice to exclude them within the current set-up of Bitcoin blockchain.

\subsubsection{Determining time investments}
\label{sec:determining_time_investments}

Let $\{\txset_l\}_{l\in\isetl},$ where $\isetl := \{0,\dots, L-1\}$, be any segmentation of mempool of Bitcoin transactions by Algorithm \ref{alg:tx_segmentation}. By $\block_l$ we depict the sub-blocks of Figure \ref{fig:hbs_svc_block_scheme}, such that  $\block_l\subset\txset_l.$ All transactions within a segment $\txset_l$ enjoy the same miner effort if included in a block $\block_l$, i.e. $\ti(t) = \ti_l$ for all $t\in \txset_l$ and $l \in \isetl.$ Using \refe{eq:def_security} and \refe{def:tx_diff} we obtain:

\begin{equation}
\overline{\secc(t)}|_{\txset_l} = \frac{\thr\cdot\ti_l}{|\txset_l|}\sum_{t\in\txset_l}\frac{\bsize(t) }{\ov(t)}.
\label{eq:av_security_t_l_hsb}
\end{equation}

Further, all transaction $t$ in each  $\txset_l, l \in \isetl$ have by construction similar output values per bit  $\vpb(t) = \ov(t)/\bsize(t)$. Consequently, we may substitute in the average value $\overline{\vpb(t)}|_{\txset_l}$ of $\vpb(t)$ and we obtain

\begin{equation}
\overline{\secc(t)}|_{\txset_l} \approx  \frac{\thr\cdot\ti_l}{\overline{\vpb(t)}|_{\txset_l}},
\label{eq:average_secc}
\end{equation}
which is a simple formula containing only $\thr$ and $\ti_l$ and average values of $\secc$ and $\vpb$ per segment $\txset_l.$ 

Under Assumption \ref{ass:same_security} we strive to secure the same security regardless of transaction value per bit $\vpb,$ which translates to:

\begin{equation}
\overline{\secc(t)}|_{\txset_k} = \overline{\secc(t)}|_{\txset_l}\quad\forall\ k,l\in \isetl,
\label{eq:condition_average_sec_assum}
\end{equation}
i.e. Assumption \ref{ass:same_security} is satisfied in a per sub-block fashion. Using \ref{eq:average_secc} we obtain practically the final condition which allows us to determine the time investments $\ti_l, l\in\isetl$ of the individual sub-blocks:
\begin{equation}
\frac{\thr\cdot\ti_k}{\overline{\vpb(t)}|_{\txset_k}} = \frac{\thr\cdot\ti_l}{\overline{\vpb(t)}|_{\txset_l}}\quad\forall\ k,l\in \isetl,
\label{eq:condition_average_sec_assum_ti}
\end{equation}
i.e. we look for a unknown $\cti\in \RR^+$ such that 
\begin{equation}
\ti_l \left[\frac{\mbox{s}}{\mbox{b}}\right] = \cti \left[\frac{\mbox{s}}{\bitcoinA}\right] \cdot \overline{\vpb(t)}|_{\txset_l} \left[\frac{\bitcoinA}{\mbox{b}}\right]\quad\forall l\in \isetl,
\label{eq:condition_average_sec_assum_final}
\end{equation}
where $\cti$ plays the role of $c$ from the last lines of Section \ref{sec:transaction_security} under the specific circumstances of this section. It belongs to the only one factor of difficulty $\diff$ we control here: the time investment $\ti.$ It is common to all levels $l\in\isetl$ and its unit (i.e. $s/\bitcoinA$) helps us to understand that it essentially represents computational time invested to secure one $\bitcoinA$ of transacted wealth.

Assume for a moment that a transition from a current block to a new multi-block is possible, so to say, at once. We can easily compute $\cti$  from current global difficulty $\Diff$. First, based on the definition of transaction difficulty \refe{def:tx_diff}, we spread the current hashrate in the correct way among the sub-blocks $\{\block_l\}_{l\in \isetl}$ from Figure \ref{fig:hbs_svc_block_scheme} while achieving the same expected block creation time of $600$s for the super-blocks by asking that:
\begin{equation}
\thr\sum_{l\in\isetl}\ti_l\sum_{t\in \block_l} \bsize(t) = 2^{32}\Diff.
\label{eq:distribution_of_hr_hsb}
\end{equation}
Thus it holds
\begin{equation}
600s \mathop{=}\limits^{\refe{eq:difficulty}} \frac{2^{32}\Diff}{\thr} = \sum_{l\in\isetl} \ti_l\cdot\bsize(\block_l) \mathop{\approx}\limits^{\refe{eq:condition_average_sec_assum_final}}
\cti\sum_{l\in\isetl} \overline{\vpb(t)}|_{\block_l} \cdot \bsize(\block_l),
\label{eq:c_eta_derivation_1}
\end{equation}
where the last relation is only approximation, since blocks $\{\block_l\}$ are chosen from mempool's segmentation $\{\txset_l\}$ by miners. It can be expected that miners will prefer the transactions paying higher fees per bit in each sub-block $\block_l$, which should be correlated to the value $\vpb(t)$ being transacted. One could redefine \refe{eq:condition_average_sec_assum_final} by using blocks instead of the mempool segmenation, but no mathematical trick can cover the fact that miners are responsible for creation of valid blocks and distributions of transactions in $\{\block_l\}_{l\in\isetl}$ can potentially differ from those in mempool segmentation $\{\txset_l\}_{l\in\isetl}$.  

Assume further, that $\overline{\vpb(t)}|_{\block_l}$ can be approximated as
\begin{equation}
\overline{\vpb(t)}|_{\block_l} \approx \frac{\overline{\ov(t)}|_{\block_l}}{\overline{\bsize(t)}|_{\block_l}}
\label{eq:approx_average_beta}
\end{equation}
with a sufficient accuracy. Then we can continue simplifying \refe{eq:c_eta_derivation_1} as follows:
\begin{equation}
\cti\sum_{l\in\isetl} \overline{\vpb(t)}|_{\block_l} \cdot \bsize(\block_l) \approx
\cti\sum_{l\in\isetl} \frac{\overline{\ov(t)}|_{\block_l}}{\overline{\bsize(t)}|_{\block_l}}|_{\block_l} \cdot \overline{\bsize(t)}|_{\block_l} |\block_l| = 
\cti\sum_{l\in\isetl} \ov(\block_l),
\end{equation}
i.e. 
\begin{equation}
\cti \approx \frac{600s}{\sum_{l\in\isetl} \ov(\block_l)} = \frac{2^{32}\Diff}{\thr \sum_{l\in\isetl} \ov(\block_l)}.
\label{eq:c_eta_derivation_2}
\end{equation}

We can choose \refe{eq:c_eta_derivation_1} or \refe{eq:c_eta_derivation_2} to determine $\cti$, depending on the accuracy of \refe{eq:approx_average_beta}. The second expression explains very clearly what $\cti$ represents: the ratio of expected super-block creation time $600s$ to the total transacted value of this super-block $\ov(\{\block_l\}_{l\in\isetl})$, i.e. in short time of all miners to wealth of all users.

The relation \refe{eq:c_eta_derivation_2} also suggests how the computation of $\cti$ should be implemented in practice. Its expression contains difficulty $\Diff$, unknown and variable total hash rate of Bitcoin network $\thr$ plus again variable total transacted value $\ov(\{\block_l\}_{l\in\isetl})$. Just like in the case of $\Diff$, $\cti$ can be computed only probabilistically. There is no reason to change what is working very well. We propose adjusting $\cti$ every 2016 blocks with the goal to achieve average super-block time of 10 minutes.

Now, we can return back to relation \refe{eq:condition_average_sec_assum_final} with an already determined $\cti$. Purely theoretically, one could directly use \refe{eq:condition_average_sec_assum_final} to calculate $\ti_l$ for each new super-block. However $\overline{\vpb(t)}|_{\block_l}$ can fluctuate a lot, especially at the higher levels, see e.g. level 0 in Table \ref{table:segmented_txs}, where we have only 615 transaction in 493 blocks, i.e. sometimes level 0 may be even empty. Second, since $\ti$ controls difficulty of the sub-blocks, which are essentially normal Bitcoin blocks, all the attacks possible when adjusting difficulty per block (and potentially more) would be possible \cite{nullc2021}. Indeed, difficulty is fundamentally a security mechanism and so are all the time investments $\{\ti_l\}|_{l\in\isetl}$. Thus we again propose updating them in the same fashion as difficulty of normal blocks is updated - every 2016 blocks:
\begin{equation}
\ti_l  = \frac{\cti}{2016}\sum_{i=0}^{2015}\overline{\vpb(t)}|_{\block_{l, i}}\quad\forall l\in \isetl,
\label{eq:ti_l_works_as_dificulty}
\end{equation}
where $\block_{l, i}, i=0, \dots, 2015$ for each $l\in\isetl$ are the last $2016$ sub-locks in the blockchain.

%\subsubsection{Discussion and analysis}
%\label{sec:disc_and_analysis}

%We assume that $\ti_l \ne 0$ for all $l\in\isetl$, i.e. at least one block in %\refe{eq:condition_average_sec_assum_final} is always non-empty and contains a transaction with positive %output value. It should not be a problem to fulfill this assumption in practice. For miners it is %preferable to have non-empty sub-blocks since then they generate extra income which supplements their %block rewards. And the users will be prepared to pay for adequate security and quick inclusion into a %block. We think that due to market forces an equilibrium among fees and $l$ will be quickly found. We also %expect the equilibrium to be rather static and any changes of $\{\ti_l\}$ quite smooth since  $\txset$ is %already very rich, see Figure \ref{fig:tx_vpb}. 

Let us compute $\cti$ and $\ti_l$ for our dataset of transactions, using function \verb|compute_c_eta_and_eta| from Algorithm \ref{alg:compute_c_eta_and_eta}. The estimate for $\cti$ based on \refe{eq:c_eta_derivation_1} reads
\begin{equation}
\hat{\cti} \approx 0.036 \frac{\mbox{s}}{\bitcoinA}
\label{eq:eta_data_estimate}
\end{equation}
and
\begin{equation}
\hat\ti \approx (\num{0.13}, \num{8.2e-04}, \num{2.3e-05}, \num{4.7e-07}, \num{1.6e-08}, \num{2.6e-10})
\label{eq:ti_l_data_estimate}
\end{equation}
with the values representing time in seconds spent by all the miners to secure one bit of information. These values are corresponding to levels from Table \ref{table:segmented_txs}. When we compare them with the corresponding estimated value of time investent for current Bitcoin blocks $\hat\ti_B \approx \num{7.15e-5}$ using \refe{eq:eta_B_estimate}, we see that the transactions in the first 2 sub-blocks would enjoy higher immediate security guaranties than nowadays. Looking at Table \ref{table:segmented_txs}, it concerns transactions starting from $\ov = 6.8\bitcoinA.$ On the other hand, a complex transaction even moving as much as $\ov = 1600\bitcoinA$ could be included in a sub-block $l=2$ and pay fees more than 3 times lower than today, see \refe{eq:fees_per_bit_level_l}.

Based on \refe{eq:ti_l_data_estimate}, we can estimate the time required to sign the individual levels for an average block as
\begin{equation}
\hat t_l \approx \hat{\ti_l}\cdot \overline{\bsize(\block_l)}
\label{eq:t_l_data_estimate}
\end{equation}
and thus using function \verb|compute_time_per_level| from Algorithm \ref{alg:compute_c_eta_and_eta} we get:
\begin{equation}
\hat t \approx (429, 124, 44.5, 2.85, 0.0247, \num{5.89e-06}).
\label{eq:time_per_level_L_6}
\end{equation}

Here we have to, for the first time, distinguish two possibilities: whether miners are expected to broadcast each sub-block $\block_l, l\in\isetl$ or they should work on sub-blocks independently from each other and broadcast only the final multi-block. %Choosing between the options has an impact on the time of transaction confirmation. 
When miners broadcast the sub-blocks, transactions are verifiable sooner than when they only broadcast the finalized multi-blocks. However since the difference between the mining times of two successive sub-blocks belonging to a same level $l\in\isetl$ is on average $10$ minutes, the confirmation will practically never be immediate and the upper bound is on average 10 minutes even if the next sub-block $\block_l$ is able to contain all waiting $l-$level transactions. From the usability point of view, there is thus no significant difference between these options. Since in this section we consider that all the transactions are seen by all the miners (and full nodes), one might better opt to broadcast only the finalized blocks. If we were able to transition from a current block to a multi-block at once, this would allow us to preserve the current broadcasting implementation and to introduce a high number of levels $L$.

However, to be able to scale PoW in the next sections, starting from Section \ref{sec:hbs_binary_tree}, we have to broadcast the individual sub-locks. One observes in \refe{eq:time_per_level_L_6} that for the current distribution of transaction output values $\ov$, more then 4 levels of sub-blocks are practically unthinkable, since the last sub-blocks would be created in fractions of a second and broadcasting them safely takes certainly much longer. Technically, based on our understanding of state of the art, three levels are certainly possible. E.g. Ethereum's blockchain has 15 seconds block time. The fourth level might be a challenge. One may consider a hybrid approach such that the last sub-blocks are mined in a sequence and only then broadcasted together, here e.g. sub-blocks $\block_l, l\in\{2,3,4,5\}.$ This is very easy to implement, namely miners broadcast only when sufficient work has been done.
%We intrust the fine treatment of technical limitations to the Bitcoin community.

But why are the mining times for the last sub-blocks in \refe{eq:time_per_level_L_6} so small? The transition from a single block to a multi-level block is a classical example of the chicken and egg problem. Since transactions $t$ with small $\ov(t)$ on-chain are expensive, they are rare. 

\subsubsection{Value transfer downstream}
\label{sec:value_transfer_downstream}
Whatever form the broadcasting scheme takes, we propose the following iterative approach to gradually increase the number of possible levels $L$. First, introduce a viable number of initial levels - e.g. three. As argued in Section \ref{sec:fees}, users should pay fees proportional to difficulty $\diff(t),$
which in settings of this section equals
\[\diff(t) = \thr\cdot\ti_l\cdot\bsize(t)\]
for a transaction $t$ from a block $\block_l, l\in\isetl.$ It is, the \emph{fees per bit} $\tfpb_l$ at level $l\in\isetl$ should be proportional to $\ti_l$:
\begin{equation}
\tfpb_l \left[\frac{\bitcoinA}{b}\right] \propto \ti_l, 
\label{eq:fees_per_bit_level_l}
\end{equation}
which under validity of Assumption \ref{ass:same_security} is also proportional to the average value per bit transacted within the block $\block_l$, see \refe{eq:condition_average_sec_assum_final}. 

Looking at \refe{eq:ti_l_data_estimate}, transactions of a second sub-block $\block_1$ should thus pay $0.13/\num{8.2e-04} \approx 159$-times less $\tfpb$ and those within a third sub-block $\block_2$ approximately $5650$-times less than transactions in a first sub-block $\block_0$, whatever the corresponding absolute values would be. These relative differences represent strong economic incentives to choose a lower level blocks over a higher level one. Moreover, if accepted by miners, a high value transaction belonging normally e.g. in a highest $0$-level block can be included in a lower lever block and user can choose to wait for a confirmation with a number of $0-$level blocks. These two forces would certainly move a lot more wealth from higher levels to lower levels, filling those in and consequently allowing to increase the number of levels $L.$ A mechanism for automatic adaptation of $L$ with some save maximum of e.g. $16$\footnote{From 21000000 $\bitcoinA$ to 1 sat in $\lg$-scale.} can be included into protocol. In essence, one could e.g. ask that the last sub-block time is at least a predefined constant $t_{\min}(\mathrm{block})$:
\begin{equation}
t_{L-1} \approx \ti_{L-1}\cdot \overline{\bsize(\block_{L-1})}> t_{\min}(\mathrm{block}).
\label{eq:last_level_L_condition}
\end{equation}

\subsubsection{Block rewards}
\label{sec:block_rewards_hbs}
If sub-blocks $\block_l, l\in\isetl$ are never broadcasted, only the final multi-block, no changes to reward scheme are necessary. The new multi-block takes the place of the current Bitcoin block.

If sub-blocks are broadcasted, it is necessary to distribute block rewards $\br$ among the individual levels in a multi-block from Figure \ref{fig:hbs_svc_block_scheme}. First, the expectation that at most 21M of coins will ever exist should be preserved. To comply, it is sufficient to apply the current block reward scheme at the multi-block level, since we aim at the same average multi-block creation time of $600$s . This reward then should be split among the individual sub-block based on their average creation time. We propose assigning rewards 
\begin{equation}
\br_l  := \frac{\hat t_l}{\sum_{l\in\isetl} \hat t_l}\br
\label{eq:reward_at_level_l}
\end{equation} 
for each $l\in\isetl,$ where $\hat t_l$ are estimators of average time to mine sub-block $\block_l$, see \refe{eq:t_l_data_estimate}, computed at the same time and in the same way as time investments $\ti_l, l\in\isetl$ are computed in \refe{eq:ti_l_works_as_dificulty}, i.e.

\begin{equation}
\hat t_l := n_l\cdot\overline{\bsize(\block_l)} \approx \frac{\ti_l}{2016}\sum_{i=0}^{2015}\bsize(\block_{l, i})\quad\forall l\in \isetl,
\label{eq:estimator_average_block_time}
\end{equation}
where again $\block_{l, i}, i=0, \dots, 2015$ for each $l\in\isetl$ are the last $2016$ sub-locks in the blockchain.

\subsubsection{How to transition from block to multi-block?}
\label{sec:transition}
When switching from the current Bitcoin block structure to the multi-layered solution of Figure \ref{fig:hbs_svc_block_scheme}, continuous operation has to be ensured. This is not a trivial task and it has to be carefully planned and its full treatment is outside of the scope of this publication. Nevertheless, in our view, it should be done in a maximally not intrusive way where new functionality is a possibility, not a dictate. Here, we present such an approach on conceptual level.

Bitcoin represents among other important things a breakthrough innovation, store of accumulated trust and even a promise of a better future for mankind. Naturally, its development is less agile and more conservative and careful than those of other relatively never blockchains. For this reason, what we propose in this paper is a stepwise and systematic road map of metamorphosis, which makes all the three main parties, i.e, users, miners and developers, comfortable enough to undergo the transition, in our view a vital one.

\begin{figure}[t]
	\begin{tikzpicture}
	
	\node[draw=none] at (0.5, \hsbFigOneShift)    (dots_1) {$\dots$};
	\node[draw] at (2, \hsbFigOneShift)	(block_1) {block};
	\node[draw] at (4, \hsbFigOneShift)	(block_2) {multi-block};
	\node[draw] at (6, \hsbFigOneShift)	(block_3) {block};
	\node[draw] at (8, \hsbFigOneShift)    (block_4) {multi-block};
	\node[draw] at (10, \hsbFigOneShift)	(block_5) {block};
	\node[draw=none] at (11.5, \hsbFigOneShift)(dots_2){$\dots$};
	\draw[-stealth] (dots_1) -- (block_1);
	\draw[-stealth] (block_1) -- (block_2);
	\draw[-stealth] (block_2) -- (block_3);
	\draw[-stealth] (block_3) -- (block_4);
	\draw[-stealth] (block_4) -- (block_5);
	\draw[-stealth] (block_5) -- (dots_2);
	
	\end{tikzpicture}
	\caption{New hybrid blockchain with multi level blocks included between the normal Bitcoin blocks.}
	\label{fig:transition_from_block_to_multiblock}
\end{figure}
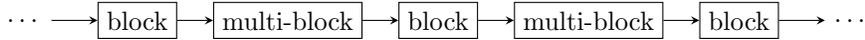

Consequently, we propose working in an additive fashion. The first elemental idea is presented in Figure \ref{fig:transition_from_block_to_multiblock}. New multi-blocks are included between the current Bitcoin blocks. Initially, the current system remains practically intact. The new multi-blocks form a small seed which will compete in usefulness with the old implementation. We are persuaded that if given chance, it will grow in time and overtake the old blocks due to its superior economic and from Section \ref{sec:hbs_binary_tree} also technical properties.

Mathematically speaking, we are constructing a homotopy from an old system to a new one, i.e. instead of the state transition function \refe{eq:transition_function} we will consider the following one:
\begin{equation}
\transf_{\lambda} = (1-\lambda)\transf_o + \lambda\transf_n,
\label{eq:transition_function_lambda}
\end{equation}    
where $\transf_o$ is the current state transition function, $\transf_n$ is the new one representing the hypothetical Bitcoin blockchain with multi-blocks only and $\lambda$ is a real parameter from interval $(0, 1)$ which ``measures'' to which degree the old blocks and to which degree the new multi-blocks are used for transactions. This homotopy is the second elemental idea.

In the beginning, only old blocks are used, i.e. $\lambda = 0.$ If new multi-blocks succeed and all transactions employ them, $\lambda = 1.$ Here, parameter $\lambda$ is controlled predominately by users and then by miners since they decide on the inclusion of transactions into blocks. Thus, if multi-blocks did not represent a winning formula for users and miners and indeed also developers in the beginning of the road, they could never succeed.

The reader can immediately see that this homotopy mechanism can be used to relatively safely and democratically introduce any new functionality. A disadvantage in our case is that the Bitcoin software has to support two or more block versions, at least temporarily. If $\lambda$ reaches almost $1$, users might be theoretically forced to transfer values to new blocks within a transition period of e.g. few years. Then support for old blocks could be deprecated. However, this would be a rather harsh measure. Consequently, the system should most probably remain backward compatible. Also, once even a single multi-block is non-empty, it caries some information about the Bitcoin state $\state$ and the state transition function $\transf_n$ can not be easily disposed of. Again, values would have to be first transfered to old blocks. In short, Equation \refe{eq:transition_function_lambda} is certainly no tool for experimentation and only worthy candidates $\transf_n$ should be considered.

Assume that one old block and multi-block pair is created on average every 10 minutes. Then
\begin{equation}
\lambda = \frac{t_n}{600s}, 
\label{eq:lambda}
\end{equation}
where $t_n$ is average time necessary to mine new multi-blocks. In practice, even if initially $t_n$ is almost zero, since the first multi-block contains only header and no transactions, one has to reserve a positive time for the block propagation. Based on Figure 4 of a rather dated paper \cite{Decker2013}, we have to expect at least $1$ms/B for propagation of en empty Bitcoin block to 90 percentile of the network, i.e approximately 80ms for a header of 80B and potentially much more.

Let us apply Algorithm \ref{alg:tx_segmentation} to obtain the corresponding segmentation for $L=2.$ What would be $t_n$ if the transactions for $l=0$ stayed in old blocks and those from $l=1$ moved to new multi-block, here with one level only? The results are presented in Table \ref{table:segmented_txs_L_2}. In comparison to Table \ref{table:segmented_txs} we added a $\overline{\bsize(t)}$-row to illustrate that small transaction are actually more ``rich''. And we added $\ti_l$- and  $t_l$-rows which are computed according to Algorithm \ref{alg:compute_c_eta_and_eta}. Otherwise, the table is generated in the same way as Table \ref{table:segmented_txs}.

\def\arraystretch{1.1}
\begin{table}
	\begin{tabular}{|c|c|c|}
		\hline
		$l$ & 0 & 1\\ \hhline{|=|=|=|}
		$|\mathcal{T}_l|$ & 245679 & 703933 \\ \hline
		$\overline{\bsize(\mathcal{T}_l)}$ & 5.33e+02 & 6.66e+02 \\ \hline
		min($\beta(t)$) & 6.39e+03 & 0.0197 \\ \hline
		max($\beta(t)$) & 2.08e+09 & 6.39e+03 \\ \hline
		$\overline{\beta(t)}$ & 1.2e+06 & 1.07e+03 \\ \hline
		min($v(t)$) & 9.65e+06 & 3.06e+02 \\ \hline
		max($v(t)$) & 6e+12 & 7.4e+09 \\ \hline
		$\overline{v(t)}$ & 3.04e+09 & 5.61e+06 \\ \hline
		$\sum v(t)$ & 7.46e+14 & 3.95e+12 \\ \hline
		$\ti_l$ & 0.000281 & 2.52e-07 \\ \hline
		$t_l$ & 5.98e+02 & 1.92 \\ \hline    
	\end{tabular}
	\caption{Distribution of transaction values after segmentation: $L=2$}
	\label{table:segmented_txs_L_2}
\end{table}

Even if the second segment $\txset_1$ of ``the mempool'' $\txset$ contains more transactions and is bigger, under the validity of Assumption \ref{ass:same_security}, the transactions included within it enjoy the same security by using the hash rate $\thr$ of the Bitcoin network for less then 2 seconds out of 600s, particularly $t_n \approx 1.92s$. And they should pay fees approximately $\ti_0/\ti_1 \approx 1120$ times lower on average then those transactions in $\txset_0.$ Maybe more importantly, the energy usage per transaction is also approximately $1120$-times lower. It means, that a small investment of miners time, i.e. approximately $1.92/598=.32\%$, represents the capacity to sign more than two thirds of all the transactions while still providing an adequate security and decreasing energy consumption of small value transactions three fold. We hope, that the huge potential effect on the ability of the Bitcoin network to serve small regular transactions on-chain is obvious. 

Is it necessary to motivate somehow the use of new multi-blocks? The short answers is no. First, it is reasonable to assume that transaction fees and environmental concerns on its own would be sufficient motivation to employ multi-block for small casual transactions. With growing confidence, higher and higher values would be transacted. If found desirable, the natural economic incentives which promote value transfer  described in Section \ref{sec:value_transfer_downstream} could be strengthened. E.g. those users using multi-blocks could enjoy a higher security $t_n + \epsilon$ for some $\epsilon>0$. However, since transaction fees are determined by market forces, such measures would be probably ineffective.

Whether an equilibrium with a $\lambda < 1$ would be found or the current blocks would stopped to be used altogether, i.e. $\lambda$ would eventually reach $1$, is a difficult question to answer. The answer depends on too many unknowns, among them on implementation details. We may only speculate that eventually the current blocks would probably survive and work as a bank for high value transactions, since \refe{eq:def_security} is only a model and in reality the relation between $\diff(t)$ and $\ov(t)$ is very likely non-linear. E.g. those making high-value transactions are maybe prepared to accept a higher risk or are only making transactions between two wallets they own and thus double-sending is of no concern or they are maybe able to wait longer for confirmation. On the other hand miners would realize this and they would include also transaction paying less than a ``fair'' linear-model based value into standard blocks.\footnote{The same reasoning applies to relations among sub-blocks of a multi-block.} But most importantly sub-blocks are essentially current Bitcoin blocks and one of them is the highest in hierarchy. Thus the standard block preceding a multi-block can be seen as the highest level sub-block for $l=-1$ and $\lambda$ is then merely a degree-of-freedom parameter which has to be found to balance the two state transition functions: $\transf_o$ and  $\transf_n$.

To summarize, we propose:

\begin{itemize} 
\item Introducing multi-blocks between the standard blocks as depicted in Figure \ref{fig:transition_from_block_to_multiblock}.
\item Adapting the consensus algorithm to accommodate the state transition function \refe{eq:transition_function_lambda}, where $\lambda$ is defined by \refe{eq:lambda}, where $t_n$ could be computed at the same time and in the same way as $\ti_l,l\in\isetl$ are computed in \refe{eq:ti_l_works_as_dificulty}, i.e. by averaging the multi-block time through last 2016 blocks.
\item Propagating the whole multi-blocks only which allows us to introduce from the beginning practically any number of levels $L$, where $L\approx16$ seems to offer good granularity.
\item Dividing the block reward $r$ honestly between $\transf_o$ and $\transf_n,$ i.e. $\lambda\br$ goes to miner finding the next multi-block and $(1-\lambda)\br$ to miner finding the next standard Bitcoin block.
\end{itemize}

The l.h.s. of Equation \refe{eq:c_eta_derivation_1} contains $t_n$ instead of $600s$. All the formulas to compute $\cti$ and $\ti$ etc. are adapted accordingly.

\subsection{Sharding levels}
\label{sec:hbs_binary_tree}
In Section \ref{sec:hbs_svc} we introduced our first hierarchical block structure (HBS), depicted in Figure \ref{fig:hbs_svc_block_scheme}. It employed the concept of extended transaction $t$, see Definition \ref{def:extended_transaction}, but only those signed by all the miners and seen by all the full nodes, i.e. $\fh(t)=1$ which implies $h(t)= \thr.$ We showed that by varying time investment $\ti$, see \refe{eq:hash_time}, this HBS allows us to successfully align three things: transaction security, fees and electric energy consumption. This was all possible because we have recognized and already partially employed the fact that all transactions are not equal. Namely, that some of them move fortunes and some of them only pennies. But since the constraint $h(t)= \thr$ remained untouched, so did the scalability, as it is succinctly expressed in Section \ref{sec:introduction}. 

Assume for a moment that the constraint $h(t)= \thr$ remained untouched and we simply allowed the new multi-blocks from Figure \ref{fig:transition_from_block_to_multiblock} to have flexible size. They could indeed accommodate more transactions. But to propagate big blocks costs time  \cite{Decker2013} and the bigger their size, the more difficult it becomes to run a full node, promoting the type of centralization we really want to avoid. 

As discussed in Appendix \ref{app:centralization_of_mining}, it is not the centralization of mining but that of control which is problematic. Miners and developers, even if they are indeed privileged, merely serve the user base of the Bitcoin community. And users cast they votes every day by buying or selling coins. This type of control is obviously fully independent on blockchain implementation including its consensus algorithm. But the second voting mechanism, the choice to run or not a full node or its specific software version has to be preserved. Altogether, simply increasing multi-blocks size is not the way to scale.

What we however can do is to distribute hashing capacity and to adapt the concept of a full node. The underlying motivation is simple. Again, when one does high-stake transactions, he will try to avoid all risks. Let us take as an example a private property transaction. Buying e.g. a house without professional legal support is probably unthinkable for most people. But when paying for ice-cream we sometimes not even count change. The difference in due diligence simply stems from different values of those transactions. Which aspect matters to us the most here is how transactions are validated and recorded. Private property ownership together with the supporting documents are usually recorded in at least a county register, i.e. in a centralized database, sometimes even a publicly searchable one. On the other hand who knows a lot about your last ice cream transaction paid in cash?  

This observation implies for Bitcoin transactions $t$ that we should require those with the highest value per bit $\vpb(t)$ to be signed by all miners and validated and recorded by everybody, i.e. all nodes. But those with lower values can be mined by only a fraction of the global hashrate $\thr$ and validated and recorded by only some fraction of the ``full'' nodes. 

We can build directly on the work presented in Section \ref{sec:hbs_svc}, where $\vpb(t)$ decreases as the level index $l\in\isetl$ increases, i.e., in effect, also the required due diligence decreases. Thus we can exploit the first degree of freedom of an extended transaction $t$, see Definition \ref{def:extended_transaction}, i.e. $\fh$, and split the individual levels in Figure \ref{fig:hbs_svc_block_scheme}. More accurately, we split the mining work, it is the transactions waiting in the mempool for confirmation among the mining nodes in the network. The validation work is split as well. The highest level for $l=0$ is not split at all. The transactions within the corresponding blocks are worked on by all the mining nodes and seen by all the full nodes. It behaves as a current bitcoin block. 

\def\hsbFigOneShift2{4.5}
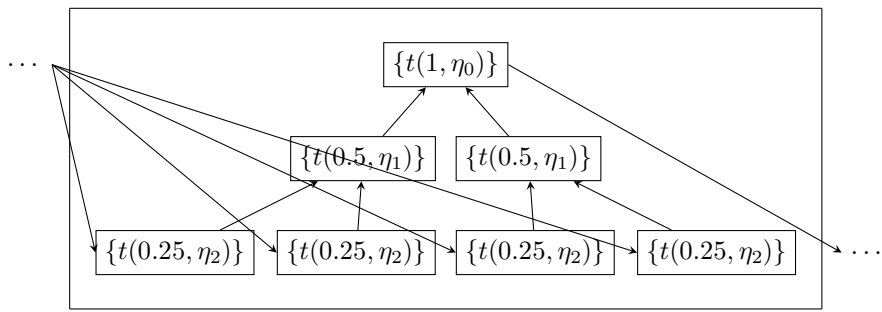
\begin{figure}[t]
	\hspace{-1.7cm}
	\begin{tikzpicture}
	\node[draw=none] at (\hsbFigOneShift2-6.5,2.5)	(level1dots) {$\phantom{\qquad\qquad\quad}\dots$};
	
	\node[draw] at (\hsbFigOneShift2,2.5)	(level1) {$\{t(1, \ti_0)\}$};
	\node[draw] at (\hsbFigOneShift2-1.1,1.25)	(level2_0) {$\{t(0.5, \ti_1)\}$};
	\node[draw] at (\hsbFigOneShift2+1.1,1.25)	(level2_1) {$\{t(0.5, \ti_1)\}$};
	\node[draw] at (\hsbFigOneShift2-3.6,0)	(level3_0) {$\{t(0.25, \ti_2)\}$};
	\node[draw] at (\hsbFigOneShift2-1.19,0)	(level3_1) {$\{t(0.25, \ti_2)\}$};
	\node[draw] at (\hsbFigOneShift2+1.19,0)	(level3_2) {$\{t(0.25, \ti_2)\}$};
	\node[draw] at (\hsbFigOneShift2+3.6,0)	(level3_3) {$\{t(0.25, \ti_2)\}$};
	\draw[-stealth] (level2_0) -- (level1);
	\draw[-stealth] (level2_1) -- (level1);
	\draw[-stealth] (level3_0) -- (level2_0);
	\draw[-stealth] (level3_1) -- (level2_0);
	\draw[-stealth] (level3_2) -- (level2_1);
	\draw[-stealth] (level3_3) -- (level2_1);
	
	\draw[-stealth] (level1dots.east) -- (level3_0.west);
	\draw[-stealth] (level1dots.east) -- (level3_1.west);
	\draw[-stealth] (level1dots.east) -- (level3_2.west);
	\draw[-stealth] (level1dots.east) -- (level3_3.west);
	
	\node[draw=none] at (\hsbFigOneShift2+6.5,0)	(levelldots) {$\dots$\phantom{\qquad\qquad\quad}};
	
	\draw[-stealth] (level1.east) -- (levelldots.west);
	
	\draw (\hsbFigOneShift2-5,-0.75) -- (\hsbFigOneShift2-5,3.25);
	\draw (\hsbFigOneShift2-5,3.25) -- (\hsbFigOneShift2+5,3.25);
	\draw (\hsbFigOneShift2+5,3.25) -- (\hsbFigOneShift2+5, -0.75);
	\draw (\hsbFigOneShift2+5,-0.75) -- (\hsbFigOneShift2-5, -0.75);
	\end{tikzpicture}
	\caption{Sharding of individual levels for $L=3$.}
	\label{fig:hbs_tree}
\end{figure}

One possible realization of the corresponding HBS is depicted for $L=3$ in Figure \ref{fig:hbs_tree}. Here, we have employed a binary tree structure but any number of children per a non-leaf node $\noc$ is possible. For simplicity but without any loss of generality we assume that $c=2$ except when otherwise noted.

Formally, we still map the transactions from $\txset$ based on their $\vpb$ and/or based on what users desire to one-dimensional space with the main parameter being $l\in\isetl$. The individual sub-blocks $\block_{l, s}, s\in\{0,\dots, 2^l-1\}$ for a certain level $l\in\isetl$ are principally equivalent to each other with respect to security, since they employ the same time investment $\ti_l.$

Even if the idea is simple, the set up is much more complex than that of Section \ref{sec:hbs_svc}. And to reach any form of consensus in the community will be probably an equally complex and challenging process. Nevertheless, we again suggest introducing any new functionality using the homotopy mechanism from Section \ref{sec:transition}, which makes the transition process a rather safe one.  It is even possible that the development of this section and Section \ref{sec:hbs_svc} can be aggregated to avoid a two stage transition process.

%What follows is a conceptual presentation of a road map to exploit the first degree of freedom of

\subsubsection{Distribution of hashing resources and its security}

First, in the HBS from Figure \ref{fig:hbs_tree}, we assume implicitly that the hashrate $\thr$ can be uniformly divided among the shards. This can be a problem especially for a higher number of levels $L$, when some of the biggest miners can have capacity higher then is necessary to sign a single shard $\block_{l,s}$. We expect that such sophisticated miners have no problem to split their mining power among the individual shards if necessary.

Moreover, nowdays hashing capacity is mostly concentrated into mining pools. These again sophisticated actors would have to accommodate their business model and to distribute the hashing power according to the new multi-layer sharded block structure.

We propose assigning miners to shards $\block_{l, s}, s\in\{0,\dots, 2^l-1\}$ based on a hash of their IP address or any other similar unique network identifier. One can take e.g. few leading bits/bytes of this hash to assign them a shard to work on. For convenience, we define an abstract shard function $\shardf$ which does precisely this:
\begin{equation}
\shardf(l, \mathrm{hash}(IP)): (l, IP) \to s.
\label{eq:shard_function}
\end{equation}

A Python implementation of \refe{eq:shard_function} is presented in Algorithm \ref{alg:shard_func}. A miner obtains not only a shard $s_l$ at level $l$ but a branch $(s_0, s_1, s_2, \dots, s_l), s_0 = 0$ of the tree HBS, he should be working on or better to which he belongs. This branch is such that each $\block_{i, s_i}$ is a parent sub-block of shard $\block_{i+1, s_{i+1}}$ for $i\in\{0,\dots,l-1\}.$ The probabilities to belong to a certain shard at any level $i\in\{0,\dots,l\}$ are uniformly distributed.

Even if an attacker achieves more than $50\%$ hashrate for a particular shard $\block_{l,s}$ for a sufficiently high level $l \in \isetl$, his ability to double spent is very limited. After broadcasting its block, miners from a higher level $l-1$ are already including the hash of its block into their block secured by a much higher accumulative hashrate. To counterfeit any block later is also impossible, since after broadcasting it is available to all miners from level $l-1$ and also full nodes keeping records of this part of the tree HBS. In the end, the whole multi-block from Figure \ref{fig:hbs_tree} is secured by cumulative hashrate of the whole network $\thr$ when the corresponding root block $\block_{0,0}$ is mined. Consequently, the whole tree would be at least as secure as the current Bitcoin implementation.

One can also consider an attacker having enough hashrate to theoretically double spend at a certain level $l$, but to attack at a lower level $l+1.$ Such attacks are in our view highly unlikely since, simply put, monetary rewards are exponentially decreasing function of $l$ and one time only at level $l+1$. This is in sharp contrast to continuously flowing block rewards and transaction fees at the level $l$ at which the miner normally operates, together with a possible crushing affect on the BTC market price.

To further strengthen security, the function \refe{eq:shard_function} may be randomized. The hash can include a random string determined e.g. during assemble of a previous multi-block by all the miners in a distributed fashion, so that the shards assigned by \refe{eq:shard_function} are always different and not fixed. The following algorithm to generate such a global nonce comes to mind:

\begin{algorithm}[H]
	\caption{An algorithm for generation of a global random nonce.}
	\label{alg:tx_global_random_number}
	\DontPrintSemicolon
	\SetAlgoLined
	\KwData{$L\in\NN^+$; 
%set $\{R_{l,s}| l\in\isetl; s\in\{0,\dots, 2^l-1\}\}$ of random numbers stored in the last set of blocks}
	}
	\KwResult{A random nonce $R_{0,0}$ saved at the root block $\block_{0,0}$.}
	\While(\tcp*[f]{an abstract loop corresponding to multi-blocks}){1}{
		\For{$l\gets L-1$ \KwTo $0$}{
			\For(\tcp*[f]{non-blocking}){$s\gets 0$ \KwTo $2^l-1$}{
				\If{miner should mine $B_{l,s}$}{
					$\cdots$\;
					\If{$l = L-1$}{
						$R_{l,s}$ = hash(rand())\;	
					}
					\If{$l < L-1$}{
						broadcasting and validating blocks $\block_{l+1, 2s}$, $\block_{l+1, 2s+1}$\;
						$R_{l,s}$ = hash(rand(), $R_{l+1, 2s}$, $R_{l+1, 2s+1}$)\;
					}
					Add $R_{l,s}$ to the header of $B_{l,s}$\;
				    mine $B_{l,s}$\;
				    $\cdots$\;
				}
			}
		}
	}
\end{algorithm}

The generation of this random nonce is local as each miner of a block $\block_{l,s}$ for some $0\le l<L-1$, $s\in{0, \dots, 2^l-1}$ needs to know only the lower level nonces from two children blocks $\block_{l+1, 2s}$ and $\block_{l+1, 2s+1}.$ But the result $R_{0,0}$ is available to everybody, since $\block_{0,0}$ is mined by all the miners and validated and recorded by all the full nodes. Instead of \refe{eq:shard_function} we use 
\begin{equation}
\shardf(l, \mathrm{hash}(IP, R_{0,0})): (l, IP) \to s.
\label{eq:shard_function_secured}
\end{equation}

Admittedly, a sophisticated attacker can own a range of IPs and to route traffic always through one which maps to a fixed shard $s,$ and to e.g. become a miner with more then $50\%$ hashrate for this shard. However, information about IPs is publicly available and it is quite easy to propose countermeasures. E.g. it is more costly to control static IPv4 addresses with different 2 or 3 leading bytes. I.e., to create a distribution of IPs which seems randomly chosen is difficult. Altogether, this measure substantially increases costs to perform any attack where concentration of hashing resources in a shard is necessary, however unlikely they are based on the previous reasoning in this section.

\subsubsection{How to assign transactions to individual shards?}

We again propose using a variation of function \refe{eq:shard_function}
to assign transactions $t$, i.e. work, to individual shards $s\in\{0,\dots, 2^l-1\}$ for some $l\in\isetl$, particularly 
\begin{equation}
\shardf(l, \mathrm{hash}(i(t))): (l, i(t)) \to s,
\label{eq:shard_function_mining_transactions}
\end{equation}
where input $i(t)$ of transaction $t$ is a single Unspent Transaction Output(UTXO) of some previous transaction.

The tree hierarchical block structure (HBS) from Figure \ref{fig:hbs_tree} that we consider in this section does not allow us to support classical Bitcoin transactions with multiple inputs, only at the highest level block $\block_{0,0}.$ Why? Since the whole point of sharding the sub-blocks $\block_l, l\in\isetl$ is to allow for independent signing of transactions and their validation in the individual shards $\block_{l,s}, s\in{0, \dots, 2^l-1}$, the mapping between UTXOs and transactions has to be unique. Otherwise, an attacker can double spend its UTXOs in possible multiple sub-blocks from set $\{\block_{l,s}\}.$ Namely, he could create different non-disjoint subsets $S_1$ and $S_2$ of UTXOs he controls, such that $\shardf(l, \mathrm{hash}(S_1))\neq\shardf(l, \mathrm{hash}(S_2))$ and those would be assigned to different shards $s_1, s_2\in\{0, \dots, 2^l\}.$ The UTXOs in $S_1\cap S_2$ would be then double spent.

To aggregate multiple UTXOs in one transaction is still possible within the HBS of Section \ref{sec:hbs_svc}. Or indeed, one classical Bitcoin transaction with multiple inputs can be substituted by a set of simple one input transactions transacting to the same output(s). Altogether, we pay at most a small loss of convenience for possibility to scale on-chain transaction capacity.

\subsubsection{Distributed computation on tree HBS}

Let us now adapt the important formulas from Section \ref{sec:hbs_svc} to the new HBS, since many of them require some revision. First, sub-blocks $\block_{l, s}$ from level $l\in\isetl$ for any shard $s\in\{0,\dots, 2^l-1\}$ employ hashrate 
\begin{equation}
\hr_l := \frac{\thr}{2^l}
\label{eq:hr_per_shard}
\end{equation} 
and consequently transactions within those blocks enjoy, instead of \refe{eq:av_security_t_l_hsb}, the following security:

\begin{equation}
\overline{\secc(t)}|_{\txset_l} = \frac{\thr\cdot\ti_l}{2^l|\txset_l|}\sum_{t\in\txset_l}\frac{\bsize(t) }{\ov(t)}.
\label{eq:av_security_t_l_thbs}
\end{equation}
By repeating the derivation from Section \ref{sec:determining_time_investments} we get the formula for time investments of individual levels under validity of Assumption \ref{ass:same_security}: 
\begin{equation}
\ti_l \left[\frac{\mbox{s}}{\mbox{b}}\right] = 2^l\cti \left[\frac{\mbox{s}}{\bitcoinA}\right] \cdot \overline{\vpb(t)}|_{\txset_l} \left[\frac{\bitcoinA}{\mbox{b}}\right]\quad\forall l\in \isetl,
\label{eq:time_investment_from_cti_and_average_beta_thsb}
\end{equation}
where again $\cti\in\RR^+$ and the interpretation of this variable does not change. We see that since hash power of miners is divided among the shards $\block_{l,s}$ for $l>0$, $2^l-$times more time than in Section \ref{sec:hbs_svc} has to be invested to compute a hash-signature of those blocks to achieve the same security for transactions inside them. When we compare the exponential grow $2^l$ with exponential decrease of $\overline{\vpb(t)}|_{\txset_l}$ for a sufficiently small $L$ (e.g. see Table \ref{table:segmented_txs}), we see that \refe{eq:monotonicity_of_ti} still holds. 

The corresponding match to relation \refe{eq:distribution_of_hr_hsb} is 
\begin{equation}
\thr\sum_{l\in\isetl}\frac{\ti_l}{2^l}\sum_{s=0}^{2^l-1}\sum_{t\in \block_{l,s}} \bsize(t) = 2^{32}\Diff,
\label{eq:distribution_of_hr_thsb}
\end{equation}
which leads, instead of \refe{eq:c_eta_derivation_1}, to
\begin{equation}
600s \mathop{=}\limits^{\refe{eq:difficulty}} \frac{2^{32}\Diff}{\thr} = \sum_{l\in\isetl} \frac{\ti_l}{2^l}\sum_{s=0}^{2^l-1}\bsize(\block_{l,s}) \mathop{\approx}\limits^{\refe{eq:time_investment_from_cti_and_average_beta_thsb}}
\cti\sum_{l\in\isetl} \sum_{s=0}^{2^l-1}\overline{\vpb(t)}|_{\block_{l,s}} \cdot \bsize(\block_{l,s}).
\label{eq:c_eta_derivation_1_thbs}
\end{equation}

But how do we compute $\cti$, since information about $\overline{\vpb(t)}|_{\block_{l,s}}$ and $\bsize(\block_{l,s})$ for each admissible $(l,s)$-pair is not available to every miner? We have to work  in a distributed fashion as it is already done in Algorithm \ref{alg:tx_global_random_number}. Assume, that we again want to compute $\cti$ based on the last $2016$ multi-blocks as it was done in Section \ref{sec:determining_time_investments}. Then we have
\begin{equation}
\cti \approx\frac{600s}{ \displaystyle
\sum_{l\in\isetl} \sum_{s=0}^{2^l-1}\frac{1}{2016}\sum_{i=0}^{2015}\overline{\vpb(t)}|_{\block_{l,s,i}} \cdot \bsize(\block_{l,s,i})},
\label{eq:c_eta_derivation_2_thbs}
\end{equation}
where $\block_{l,s,i}, i\in{0,\dots, 2015}$ are the last 2016 shards for any admissible pair $(l,s).$ Each miner belonging to shard $(l,s)$ should compute the average ``output value'' for his proposition block
\begin{equation}
\ov_{l,s} :=\frac{1}{2016}\sum_{i=0}^{2015}\overline{\vpb(t)}|_{\block_{l,s,i}}\cdot \bsize(\block_{l,s,i}).
\label{eq:c_eta_derivation_average}
\end{equation}
A miner one level above him will control the computation during validation of the children blocks and only then include the hash of this block into the header of its block.
%Consequently, it is impossible to use the secured shard function \refe{eq:shard_function_secured} for each $i\in{0,\dots, 2015}$, since to compute \refe{eq:c_eta_derivation_average} a $(l,s)-$miner needs to know a 2 week history of $(l,s)-$sub-blocks. Nevertheless,  \refe{eq:shard_function_secured} can still be used once every 2016 blocks after $\cti$ is updated, if desirable. 
We suggest to use recurrent relation for average
\begin{equation}
\ov_{l,s,i} = \frac{i}{i+1}\ov_{l,s, i-1} + \frac{1}{i+1} \overline{\vpb(t)}|_{\block_{l,s,i}}\cdot \bsize(\block_{l,s,i})
\label{eq:update_output_values_local}
\end{equation} 
to obtain an on-line estimate $\ov_{l,s,i}$ of $\ov_{l,s}$
for each $i\in\{0,\dots, 2015\}$ with any finite initial condition $\ov_{l,s,-1} < \infty.$ Now we are ready to present the distributed Algorithm \ref{alg:cti_distributed_algorithm_thbs} to compute $\cti$, see Appendix \ref{app:algorithms}.

\begin{algorithm}[th]
	\caption{An algorithm to compute $\cti$ in tree HBS from Figure \ref{fig:hbs_tree}.}
	\label{alg:cti_distributed_algorithm_thbs}
	\SetAlgoLined
	%\DontPrintSemicolon
	\SetKw{continue}{continue}
	\SetKw{break}{break}
	\KwData{$L\in\NN^+$; multi-block index number $i\in\NN, i = 0.$
		%set $\{R_{l,s}| l\in\isetl; s\in\{0,\dots, 2^l-1\}\}$ of random numbers stored in the last set of blocks}
	}
	\KwResult{Every 2016 blocks, a new $\cti$ will be available and recorded in $\block_{0,0}$.}
	\While(\tcp*[f]{an abstract loop corresponding to multi-blocks}){1}{
		\For{$l\gets L-1$ \KwTo $0$}{
			\For(\tcp*[f]{non-blocking}){$s\gets 0$ \KwTo $2^l-1$}{ 
				$\dots$\;
				\If{miner should mine $\block_{l,s, i}$}{
					\If{$l < L-1$}{
						
						broadcasting and validating blocks $\block_{l+1, 2s,i}$, $\block_{l+1, 2s+1,i}$ while also checking 
						$\ov_{l+1,2s,i}$ and $\ov_{l+1,2s+1,i}$\;
						
					}
					$\dots$\;	\tcp{block assemblage}
					$\ov_{l,s,i} \leftarrow$ use
					\refe{eq:update_output_values_local}\;
					$s_{l,s,i} = \ov_{l,s,i}$\;
					\If{$l < L-1$}{
						$s_{l,s,i} = s_{l,s,i} + s_{l+1,2s, i} + s_{l+1,2s+1, i}$\;
						}
					add $\ov_{l,s,i}$ and $s_{l,s,i}$ to block header\;
					\If{$i  = 2015 \wedge l = 0$}{
						$\cti = 600s/s_{l,s,i}$\;
						add $\cti$ to header\;
					}
					mine $\block_{l,s,i}$\;
					$\dots$\;
				}
			}
		}
		$i = i+1$\;
		$i = i\mod2016$\;
			
	}
 \end{algorithm}

Since information is distributed among the nodes of tree HBS, not only $\cti$ but any variable dependent on such information has to be computed in a distributed way as it is done in Algorithm \ref{alg:cti_distributed_algorithm_thbs}. The same is true for time investments $\ti_l, l\in\isetl.$ For those, instead of \refe{eq:ti_l_works_as_dificulty}
we have
\begin{equation}
\ti_l  = \frac{\cti}{2016\cdot 2^l}\sum_{s=0}^{2^l-1}\sum_{i=0}^{2015}\overline{\vpb(t)}|_{\block_{l, s, i}}\quad\forall l\in \isetl,
\label{eq:ti_l_works_as_dificulty_thbs}
\end{equation}
i.e. each $\ti_l$ is dependent on information from the all shards $s\in\{0,\dots, 2^l\}$. Those are connected only via the root node of the tree and thus $\ti_l$ can only be computed by applying a reduce to the sum of $\overline{\vpb(t)}|_{\block_{l, s, i}}$ through shards $s$ using intermediate results saved in the higher levels $\{0,\dots, l-1\}$ of the tree HBS. The average through recent multi-blocks for $i\in\{0,\dots,2015\}$ can be again facilitated via an on-line update analogical to \refe{eq:update_output_values_local}. Final $\ti_l, l\in\isetl$ will be saved in the root block $\block_{0,0}$ as a vector and updated at the same time as $\cti$ is updated in Algorithm \ref{alg:cti_distributed_algorithm_thbs}. Since root blocks are available to everybody, so are $\{\ti_l\}.$ Details of the algorithm are left on reader as an exercise.

Let us now reassess $\cti$ and $\ti_l$ for our dataset of transactions. First, even if the calculation of $\cti$ changes, the value remains the same, i.e. \refe{eq:eta_data_estimate} still holds. By comparing \refe{eq:condition_average_sec_assum_final} with \refe{eq:time_investment_from_cti_and_average_beta_thsb} we see that the updated estimate $\hat{\ti}_{ tree}$ of $\ti$ reads
\begin{equation}
\hat{\ti}_{tree,l} = 2^l \hat{\ti_l}\quad \forall l\in\isetl,
\label{eq:ti_tree_data_estimate}
\end{equation}
where $\hat{\ti}$ is the estimate from \refe{eq:ti_l_data_estimate}.

Second, since the sub-blocks $\block_l,l\in\isetl$ from Section \ref{sec:hbs_svc} are here divided to $2^l$ shards $\block_{l,s}$ for $s\in\{0,\dots,2^l-1\}$, those are thus on average $2^l$-times smaller and 
\refe{eq:t_l_data_estimate} implies that $2^l$ from \refe{eq:ti_tree_data_estimate} cancels out and time to sign a shard $\block_{l,s}$ is still given by \refe{eq:time_per_level_L_6}. In consequence, the hurdles to broadcast the shards are neither lower nor higher at the first sight then those with sub-blocks in Section \ref{sec:hbs_svc}. But since the shards are $2^l$-times smaller, it takes in principle $2^l$-times less time to broadcast them. While we ignored latency issues etc., the situation is clearly significantly better, namely exponentially so. 

\subsubsection{Reward and transaction fees}

Again, as argued in Section \ref{sec:fees}, users should pay fees proportional to transaction difficulty $\diff(t).$ Per shard $\block_{l,s}, l\in\isetl, s \in\{0,\dots, 2^l-1\}$, the available hashrate is given by \refe{eq:hr_per_shard}. On the other hand we have derived that the time investment $\ti$ under validity of Assumption \ref{ass:same_security}  satisfies \refe{eq:time_investment_from_cti_and_average_beta_thsb}. Altogether, for any $l\in\isetl,$ $2^l$ cancel out and the difficulty $\diff(t_1)$ for a transaction $t_1$ in $\block_{l,s}$ is identical to the difficulty $\diff(t_2)$ for a transaction $t_2$ in $\block_l$ from Section \ref{sec:hbs_svc}, i.e. fees should not change. Indeed, the fees per bit $\tfpb_l$ at level $l\in\isetl$ should be proportional to 
\begin{equation}
\tfpb_l \left[\frac{\bitcoinA}{b}\right] \propto \frac{\thr}{2^l}\ti_l \mathop{=}\limits^{\refe{eq:time_investment_from_cti_and_average_beta_thsb}}\cti \cdot \overline{\vpb(t)}|_{\txset_l}, 
\label{eq:fees_per_bit_level_l_thbs}
\end{equation}
which is, under validity of Assumption \ref{ass:same_security}, identical to what users should pay, according to \refe{eq:fees_per_bit_level_l}, for space in sub-blocks $\block_l,l\in\isetl$.

For distribution of rewards $\br$, we build on the approach presented in Section \ref{sec:block_rewards_hbs}. If the whole network hashrate $\thr$ is uniformly divided among $2^l$-shards for each level $l\in\isetl,$ we have that for each $s\in\{0,\dots, 2^l-1\}$ 
\begin{equation}
\br_{l,s}  := \frac{\hat t_l}{2^l\sum_{l\in\isetl} \hat t_l}\br
\label{eq:reward_at_level_l_thbs},
\end{equation} 
where $\hat t_l$ are estimators of average time to mine sub-blocks $\block_{l,s}$ given by
\begin{equation}
\hat t_l \approx \hat{\ti_l}\cdot \overline{\bsize(\block_{l,s})}|_{s\in\{0,\dots, 2^l-1\}}
\label{eq:t_l_data_estimate_thbs}
\end{equation}
and computed in the same distributed way as time investment estimators $\hat{\ti_{l}}, l\in\isetl$ are computed starting from \refe{eq:ti_l_works_as_dificulty_thbs}. Since $\hat t_l,l\in\isetl$ need to be available to everybody, they have to be saved in a root block $\block_{0,0}.$ Because the estimates $\hat{\ti_l}$ of $\ti_l$ should already be saved therein, it is equivalently possible to save estimates of $\overline{\bsize(\block_{l,s})}|_{s\in\{0,\dots, 2^l-1\}}$ for $l\in\isetl$ and $\hat t_l$ can be computed using \refe{eq:t_l_data_estimate_thbs}. Analogically to
\refe{eq:ti_l_works_as_dificulty_thbs}, we can write
\begin{equation}
\overline{\bsize(\block_{l,s})}|_{s\in\{0,\dots, 2^l-1\}} \approx \frac{1}{2016\cdot 2^l}\sum_{s=0}^{2^l-1}\sum_{i=0}^{2015}\bsize(\block_{l, s, i})\quad\forall l\in \isetl.
\label{eq:estimator_average_block_time_thbs}
\end{equation}
Again, \refe{eq:estimator_average_block_time_thbs} can be computed using on-line update analogical to \refe{eq:update_output_values_local}.

\subsubsection{Sharding validation and storage: minimal ``full'' node.}
\label{sec:sharding_validation_and_storage}

Without distribution of both transaction validation and storage, the recent developments of this section would be in essence meaningless. All aspects of the Bitcoin ecosystem have to be adapted to achieve scalability on-chain.

Our core distribution approach for validation and storage is very similar to that of mining. As $\vpb(t)$ decreases as the level index $l\in\isetl$ increases, we propose that a decreasing fraction of ``full'' nodes validates and stores transactions, proportionally to hashrate $\hr$ used to sign those transactions. Looking at Figure \ref{fig:hbs_tree}, it means, that all the nodes validate and store transaction for $l=0$, a half of the nodes for $l=1$ and in general 
\begin{equation}
\frac{1}{2^l}
\label{eq:fraction_of_full_nodes}
\end{equation}
of all the nodes validate and store transaction for level $l\in\isetl.$ This defines a notion of a \emph{minimal full node} (\mfn). But nodes may choose to support a bigger part of the tree HBS, even the whole one.

Of course, the fraction of network nodes \refe{eq:fraction_of_full_nodes} has to be large enough to ensure virtually $100\%$ availability of any record. This gives us an upper bound on number of possible levels. Currently, there are more then 14000 reachable Bitcoin nodes, which means that level $l=7$ is supported by more then $100$ nodes. Thus $L=8$ seems certainly feasible. To achieve e.g. $L=16$ with 100 nodes per a leaf shard $\block_{15, s}$ for any $s\in\{0,\dots, 2^{15}-1\}$ would require approximately $100\cdot2^{15} \approx 3.3$M modes. Currently, this number seems illusory, but we will show that running a \mfn{ } becomes easier as $L$ grows, which together with the obvious increase of three HBS blockchain usability should promote widespread use of decentralized finance.

We propose assigning parts of the HBS tree to candidate {\mfn}s again based on shard function  \refe{eq:shard_function}, see also Algorithm \ref{alg:shard_func}. The {\mfn}s should validate and store all sub-blocks mined within their branch. This ensures that probabilistically all parts of the tree HBS are uniformly covered with respect to both  validation and storage.

\begin{figure}[ht]
	\scalebox{0.85}{\input{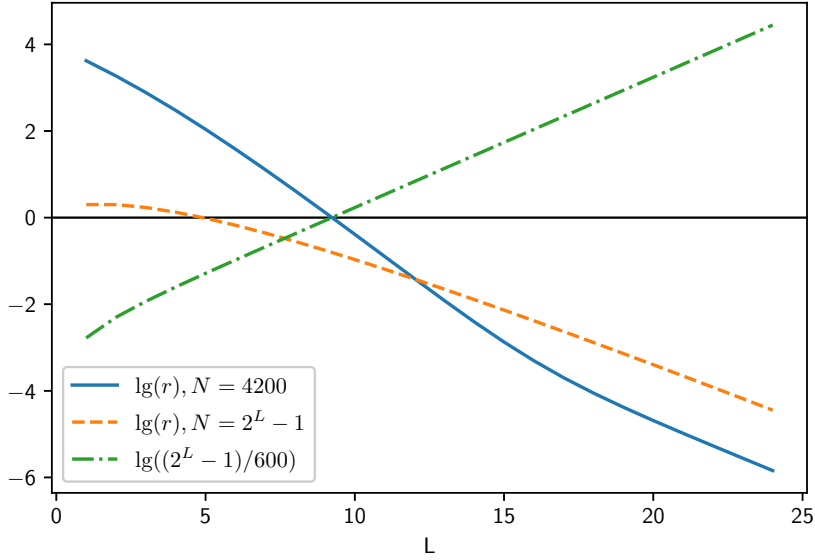}}
	\caption{Sharding efficiency: result of running Algorithm \ref{alg:sharding_efficiency}}
	\label{fig:sharding_efficiency}
\end{figure}

When a \mfn{ }has to validate transactions included in a shard $\block_{l,s}, l\in\isetl, s \in\{0,\dots, 2^l-1\}$ from the branch assigned to him by \refe{eq:shard_function}, they refer to UTXOs saved in a certain previous shard of any previous tree HBS. Thus, it has to receive from the network other sub-blocks not belonging to its branch to be able to validate those transactions. What is the number of those sub-blocks? Let us assume for simplicity that one multi-block tree contains on average $N\in\NN^+$ transactions and those are uniformly divided among the shards. Then one shard $\block_{l,s}$ for any admissible pair $(l,s)$ contains simply
\begin{equation}
\frac{N}{2^L-1}
\end{equation}
transactions, where $2^L-1$ is the number of $\block_{l,s}-$shards in the binary tree HBS with $L$ levels. To validate any of those transactions, the \mfn{ }has to download one historical branch up to some level $k<L$ from the network, where $k$ is the security level of this historical transaction. It means that it has to download, validate and store at least temporarily up to 
\begin{equation}
\frac{Nk}{2^L-1} < \frac{NL}{2^L-1}
\end{equation}
historical shards per one own shard $\block_{l,s}.$ For all $L$-shards directly assigned to a \mfn, it has to download, store and validate up to 
\begin{equation}
\frac{NL^2}{2^L-1}
\end{equation}
extra shards. Altogether, the ratio of all the shards which any \mfn{ }has to control to the alternative of storing the whole tree HBS is bounded by
\begin{equation}
r(N,L) := \frac{NL^2}{(2^L-1)^2}+\frac{L}{(2^L-1)}.
\label{eq:ratio_mfn_blocks_to_all_blocks}
\end{equation}
The scalability of \mfn{ }sharding is depicted in Figure \ref{fig:sharding_efficiency} generated by running Algorithm \ref{alg:sharding_efficiency}. The solid line represents the ratio \refe{eq:ratio_mfn_blocks_to_all_blocks} plotted as a function of $L$ for $N=4200$, which is approximately the number of $250$B legacy transactions savable in a current $1$MB Bitcoin block. For a still feasible $L=10$ the above expression already evaluates to a number smaller then 1. This corresponds to approximately $4200/(2^{10}-1) \approx 4$ transactions per a shard $\block_{l,s}$ for some admissible pairs $(l,s).$ 

Of course, the goal is to scale the on-chain capabilities of Bitcoin network, i.e. to increase $N.$ Expression \refe{eq:ratio_mfn_blocks_to_all_blocks} is minimized for $N=2^L-1$, which represents a hypothetical tree HBS where each node, i.e sub-block $\block_{l,s}, l\in\isetl,$ $s \in\{0,\dots, 2^l-1\},$ contains on average just $1$ transaction. For $N<2^L-1$ some sub-blocks $\block_{l,s}$ have to be empty. Even if this represents an intriguing possibility which should be further analyzed, we consider here only the case $N\ge 2^L-1.$ The dashed line plots the ratio \refe{eq:ratio_mfn_blocks_to_all_blocks} for $N=2^L-1$ as a function of L. We see that e.g. for already $L=15$ we achieve $r\approx0.007$, which means that a \mfn{ }has to download less than $1\%$ of total multi-block capacity. The dash-dotted line in Figure \ref{fig:sharding_efficiency} depicts the corresponding theoretical capacity of the network in number of transaction per second, under assumption of average multi-block time of $600s$.
E.g. for $L=20$ we would achieve approximately $1750$ transactions per second, which is close to the average number of VISA network transactions per second.

\subsubsection{Information routing on tree HBS}

The last aspect we shortly comment on is the question of information routing within the network. It is a rather complex problem whose full treatment lies outside of the scope of this conceptual publication. Nevertheless, the communication has to be scaled as well, otherwise it would be a bottleneck. It seems wise to employ the natural topology of tree HBS. 

Let us image e.g. that a node wants to broadcast a transaction that should be signed at a certain level $l, 0<l\le L-1$. By computing \refe{eq:shard_function_mining_transactions} for its transaction the broadcasting node identifies a branch $(s_0, s_1, \dots, s_{l})$ of the tree HBS to which this transaction belongs. From its routing table, it can identify a safe number of nodes which belongs to the corresponding shard $s_1$ based on \refe{eq:shard_function} at level $l=1.$ Those nodes should have a sufficient knowledge about their part of the tree HBS, saved in their routing table. They will then broadcast the transaction to the right nodes from shard $s_2$ at level $l=2.$ Those can broadcast it among them and route to those from $s_3$ at $l=3$ etc. It is, the communication happens along the corresponding branch of tree HBS. The broadcasting node should receive a confirmation that at least one mining node or a desirable number of them working at level $l$ is aware of the corresponding transaction. 

The routing table may or may not be distributed. E.g. we estimated in Section \ref{sec:sharding_validation_and_storage} that we need $3.3$M modes for $L=16$ at a minimum of 100 nodes per a leaf shard. A routing table would then have a size proportional to approximately $13$MB for 4B per an IP address, which is arguably not problematic. However, with e.g. billions of nodes in mind, it is still sensible to consider a non-distributed routing table? If a node knows $n$ random peers, probability that all those do not belong to a certain shard at level $l$ is given by
\begin{equation}
p(n)=\left(1 - \left(\frac12\right)^l\right)^n,
\end{equation}
which yields that it has to know at least
\begin{equation}
n=\frac{\lg (1-q)}{\lg\left(1 - \left(\frac12\right)^l\right)}
\end{equation}
nodes to be sure that with probability $q = 1-p$ at least one of $n$ nodes belongs to the considered shard. Then the broadcasting node could communicate with this node directly with a high probability, without using the briefly proposed schema above. E.g. for $p = \num{1e-10}$, $L=16$ and the lowest level $l=L-1$ one gets ``only'' $n\approx 754500.$ However for $p = \num{1e-10}$, $L=24$ and the lowest level $l=L-1$ one gets more than $193$M nodes. And worse, for unsatisfactory $p = 0.1$, $L=24$ and the lowest level $l=L-1$ one gets still more than $19$M nodes. To conclude, a distributed routing table is a must for a high number of levels $L.$ But we think that all aspects of any next evolutions of Bitcoin network should be implemented as efficiently as possible from day one. Of course, taking into account all criteria, such as security and simplicity etc. A communication scheme which exploits the tree HBS topology fully is desirable.

\subsection{Concurrent levels}
\label{sec:hbs_concurrent_levels}

In Section \ref{sec:hbs_binary_tree} we proposed a tree hierarchical block structure (HBS), depicted in Figure \ref{fig:hbs_tree}, which could theoretically allow us to scale Bitcoin blockchain to thousands of on chain-transactions per second, while building on and respecting previous developments in Section \ref{sec:hbs_svc}, where, recalling again, we successfully align three things: transaction security, fees and electric energy consumption. Moreover all processes including mining, validation and storing of blocked transactions are distributed and nodes has to handle only a small fraction of the full state $\state$.

However, at least one significant issue remains unresolved. Among any too sets of confirmed transactions recorded in sub-blocks $\block_{l,s}, s \in\{0,\dots, 2^l-1\}$ of two successive ``multi-blocks'' at any level $l\in\isetl$, there lie all sub-blocks of the current multi-block for higher levels $\{l-1, l-2, \dots, 0\}$ and the next multi-block for lower levels $\{L-1, L-2, \dots, l+1\}.$ Confirmation times of those sets of transactions differ on average by multi-block creation time $T$, which we assumed to be equal to the current block time of $600s$. Consequently, a user has to wait at most 10 minutes for its transaction to be included in a sub-block $\block_{l,s}$ even if the tree HBS from Section \ref{sec:hbs_binary_tree} has enough capacity to store all the waiting transactions. In short, transactions are not timely.

There exist few ways how to mitigate this problem. Obviously, transactions which need to be really quick can employ some off-chain mechanism such as Lightning Network \cite{Beres2020}. While we fully recognize the importance of ``perpendicular'' approaches, we would still prefer an on-chain solution, if one can be found. 

The obvious direct approach is to decrease $T$ to few seconds, e.g. $T = 5s.$ Based on behavioral research \cite{Nah2004}, tolerable waiting time for information retrieval is $2s$ and $10s$ is about the limit for keeping user's attention focused. These ``multi-block'' times are very short and push the tree HBS from Section \ref{sec:hbs_binary_tree} to uncharted waters. Surely, if the number of transactions per second is kept constant, size reduction of sub-blocks $\block_{l,s}, l\in\isetl,$ $s \in\{0,\dots, 2^l-1\},$ would be proportional to reduction in $T.$ However, it is an open question how efficiently can all the aspects of state transition \refe{eq:transition_function} be rescaled, especially with regards to information synchronization within the network. It is entirely possible that this approach can lead to an acceptable solution. A fine treatment of this possibility is however outside of scope of this publication. In the forthcoming text we propose a more systematic solution with has a far higher potential of success.
%, based on a relaxation of sequential dependence in the tree HBS from Section \ref{sec:hbs_binary_tree}.

Let us recall that Bitcoin transactions $t$ are in essence signed by expenditure of computation resources, expressed by hash rate $\hr(t)$, for a certain time, expressed by time investment $\ti(t)$ defined in \refe{eq:hash_time}. Transaction security can only be controlled by controlling those two degrees of freedom (DOFs), see Section \ref{sec:ext_tx}, particularly Equation \refe{eq:difficulty_control}.
The tree HBS from Section \ref{sec:hbs_binary_tree} depicted in Figure \ref{fig:hbs_tree} has been the first HBS we proposed in this paper which is designed to simultaneously employ those two DOFs. While it certainly does not represent a unique possibility with respect to this, it is a natural extension of the multiblock idea from Section \ref{sec:hbs_svc}. So how to achieve timeliness of transactions while keeping this tree HBS or altering it only minimally and thus preserve everything we have achieved so far?

The answer is concurrency. One has to relax the sequential dependence in Figure \ref{sec:hbs_binary_tree}. A transaction $t$ at a level $l\in\isetl$ should not wait for all sub-blocks of the current multi-block for higher levels $\{l-1, l-2, \dots, 0\}$ and of the next multi-block for lower levels $\{L-1, L-2, \dots, l+1\}$ from its branch $\{s_0, s_1, \dots, s_{L-1}\}$ to be mined. It should be included in a block $\block_{l, s_l}$ as quickly as possible. 

The idea is presented visually in Figure \ref{fig:hbs_tree_concurrent}. It depicts in the middle one blockchain containing $l-$level blocks $\block_{l, s_l, i}$ for a level $l\in\{0, \dots, L-2\}$ and $i=0,1,\dots$ surrounded by its two corresponding $l+1$ blockchains, one containing blocks $\block_{l+1, 2s_l, j}$ for $j=0,1,\dots$ and the other blocks $\block_{l+1, 2s_l + 1, k}$ for $k=0,1,\dots.$ Each of those blockchains corresponds to one node of the tree HBS depicted for $L=3$ in Figure \ref{fig:hbs_tree}, namely the nodes determined by coordinate pairs $(l,s_l)$, $(l+1, 2s_l)$ and $(l+1, 2s_l+1).$ The tree HBS with $L$ levels has together at most $2^L+1$ nodes, i.e. we obtain the same amount of concurrent blockchains. As displayed in Figure \ref{fig:hbs_tree_concurrent}, those blockchains still to a great extend respect the hierarchical topology of the underlying tree HBS from Section \ref{sec:hbs_binary_tree} and they are certainly not independent. For example, hashes of lower level blocks are included into higher level proposition blocks after performing of all necessary controls to avoid double spending etc.

What changes? We allow each blockchain $(l,s), l\in\isetl, s \in\{0,\dots, 2^l-1\}$ to grow without its blocks immediately being checked and signed by higher levels $l-1, \dots, 0.$ In essence it means, that validation by mining higher level sub-blocks/shards is postponed in time. And we also allow multiple blocks in a chain $(l,s)$ to be created without all their hashes being included in a higher level block. This does not substantially alter the tree HBS topology. A sequence of such blocks can be merely seen as a single larger block, see dashed blocks in Figure \ref{fig:hbs_tree_concurrent}. But it allows us to increase timeliness because very small blocks can be added to a chain $(l,s)$ if necessary.

\newcommand{\llb}{4.3}
\newcommand{\llpob}{3.3}
\newcommand{\yl}{1.25}
\newcommand{\ylpoa}{2.5}
\newcommand{\ylpob}{0.0}
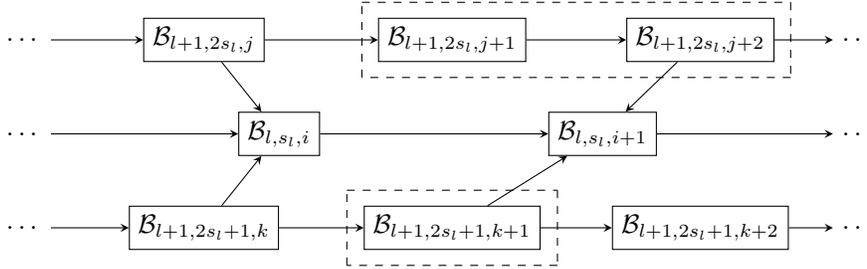
\begin{figure}[t]
	\hspace{-1.9cm}
	\begin{tikzpicture}
	\node[draw=none] at (0,\yl)	(dots_l_b) {$\phantom{\qquad\qquad\quad}\dots$};
	\node[draw] at (\llb,\yl)	(levell_0) {$\block_{l,s_l, i}$};
	\node[draw] at (2*\llb,\yl)	(levell_1) {$\block_{l,s_l, i+1}$};
	\node[draw=none] at (3*\llb,\yl)	(dots_l_e) {$\dots$\phantom{\qquad\qquad\quad}};
	
	\draw[-stealth] (dots_l_b) -- (levell_0);
	\draw[-stealth] (levell_0) -- (levell_1);
	\draw[-stealth] (levell_1) -- (dots_l_e);
	
	\node[draw=none] at (0,\ylpoa)	(dots_lp1_b) {$\phantom{\qquad\qquad\quad}\dots$};
	\node[draw] at (\llpob,\ylpoa)	(levellp1_0) {$\block_{l+1,2s_l, j}$};
	\node[draw] at (2*\llpob,\ylpoa)	(levellp1_1) {$\block_{l+1,2s_l, j+1}$};
	\node[draw] at (3*\llpob,\ylpoa)	(levellp1_2) {$\block_{l+1,2s_l, j+2}$};
	\node[draw=none] at (3*\llb,\ylpoa)	(dots_lp1_e) {$\dots$\phantom{\qquad\qquad\quad}};
		
	\draw[-stealth] (dots_lp1_b) -- (levellp1_0);
	\draw[-stealth] (levellp1_0) -- (levellp1_1);
	\draw[-stealth] (levellp1_1) -- (levellp1_2);
	\draw[-stealth] (levellp1_2) -- (dots_lp1_e);
	
	\node[draw=none] at (0,\ylpob)	(dots_lp1_t_b) {$\phantom{\qquad\qquad\quad}\dots$};
	\node[draw] at (\llpob,\ylpob)	(levellp1_t_0) {$\block_{l+1,2s_l+1, k}$};
	\node[draw] at (2*\llpob,\ylpob)	(levellp1_t_1) {$\block_{l+1,2s_l+1, k+1}$};
	\node[draw] at (3*\llpob,\ylpob)	(levellp1_t_2) {$\block_{l+1,2s_l+1, k+2}$};
	\node[draw=none] at (3*\llb,\ylpob)	(dots_lp1_t_e) {$\dots$\phantom{\qquad\qquad\quad}};
		
	\draw[-stealth] (dots_lp1_t_b) -- (levellp1_t_0);
	\draw[-stealth] (levellp1_t_0) -- (levellp1_t_1);
	\draw[-stealth] (levellp1_t_1) -- (levellp1_t_2);
	\draw[-stealth] (levellp1_t_2) -- (dots_lp1_t_e);

	% boxes	
	\def\xo{\llb+1.1}
	\def\xt{3*\llb-1.8}
	\def\yo{\ylpoa-0.5}
	\def\yt{\ylpoa+0.5}

	\draw[dashed] (\xo,\yo) -- (\xo,\yt);
	\draw[dashed] (\xo,\yt) -- (\xt,\yt);
	\draw[dashed] (\xt,\yt) -- (\xt, \yo);
	\draw[dashed] (\xt,\yo) -- (\xo, \yo);
	
	\def\xo{\llb + 0.9}
	\def\xt{2*\llb-0.6}
	\def\yo{\ylpob-0.5}
	\def\yt{\ylpob+0.5}
		
	\draw[dashed] (\xo,\yo) -- (\xo,\yt);
	\draw[dashed] (\xo,\yt) -- (\xt,\yt);
	\draw[dashed] (\xt,\yt) -- (\xt, \yo);
	\draw[dashed] (\xt,\yo) -- (\xo, \yo);
	
	%arrows among levels
	\draw[-stealth] (levellp1_0) -- (levell_0);
	\draw[-stealth] (levellp1_2) -- (levell_1);
	\draw[-stealth] (levellp1_t_0) -- (levell_0);
	\draw[-stealth] (levellp1_t_1) -- (levell_1);
	
	\end{tikzpicture}
	\caption{Concurrent chains at levels $l, l+1\in\isetl$ and some branch $s_l\in\{0, \dots, 2^l-1\}$.}
	\label{fig:hbs_tree_concurrent}
\end{figure}

Practically all results from Section \ref{sec:hbs_binary_tree} apply to this concurrent setup without a change or only with slights modifications. We lose the strict order of shards $\block_{l,s}$ within tree HBS from Figure \ref{fig:hbs_tree} but concurrent blockchains
are one-to-one mapped to individual nodes $(l,s)$ of this tree HBS, while its hierarchy is still enforced even if more loosely.  Thus this loss of strict sub-block order has only temporal implications for the finality of state $s$ when the situation is compared with the setup of Section \ref{sec:hbs_binary_tree}.

In the end, any finality of Bitcoin transactions is only probabilistic and thus prone to subjective perception, see Section \ref{sec:motivation}. Nowadays, users often wait for multiple confirmations till they consider a transaction settled and this partially irrespective of the wealth being transacted. In doing so they follow classical Byzantine fault tolerance model \cite{Nakamoto2008}. It seems that they have no other choice. At least at first glance, see \ref{app:centralization_of_mining}, this is the model governing the security of current Bitcoin blocks. All transaction within those share the same risks and guaranties. If one wants to double spend even a single transaction, he has to implicitly attack all transactions in the corresponding block. This is completely different from situation in HBSs starting already from Section \ref{sec:hbs_svc}, supported by our simple incentive based security model from Section \ref{sec:transaction_security}. Users have a choice to pick a level of security that they find appropriate based mainly on wealth being transacted. Naturally, the notion of finality is fluid as well. 

If one pays for an ice-cream or receives the payment for icecream, he is probably satisfied with the transaction being included in a low lever block for a certain low enough $l\le l_{\max},$ where $l_{\max}$ is the minimal security level both transacting parties can initially accept. They would not wait till this transaction is confirmed by creation of a highest level block $\block_{0,0}.$ E.g. the seller would accept the very small risk that finally a longer chain may be found, not referring to their transaction. As he currently accepts e.g. the risk that a very small percentage of coins are counterfeits.

On the other hand, if one receives a payment e.g. in a property transaction, he would require a much lower $l_{\max}$, such that at least a $2^{-l_{\max}}-$ fraction of nodes validate and store the data and he would consider the transaction settled only when for 
the corresponding block $\block_{l,s}$ there exists a hashpath up to the root level $l=0$ within the graph obtained by repeating Figure \ref{fig:hbs_tree_concurrent} for all pairs $(l,s), l\in\isetl, s\in\{0, \dots, 2^l-1\}.$ And even wait till a higher number of root level blocks $\block_{0,0}$ is added to the blockchain $(0,0),$ similarly to the current practice.

\subsubsection{Running a MFN node}

While for an $L$ level tree HBS, there are at most $2^L+1$ concurrent blockchains, each minimal full node (MFN), defined in Section \ref{sec:sharding_validation_and_storage}, directly handles, i.e. validates and stores, only $L$ of them, those contained in its branch of the tree HBS $\{s_0, s_1, \dots, s_{L-1}\}$, determined by \refe{eq:shard_function}. For a MFN the problem is linear in $L$. The ratio of transaction information necessary to fulfill its duties with respect to all information contained in all $2^L+1$ concurrent blockchains is estimated in \refe{eq:ratio_mfn_blocks_to_all_blocks}.

Today CPUs have mostly at least 4 cores, many supporting hyper-threading, meaning that up to $L=8$ each blockchain can be theoretically easily mapped onto a single hardware supported thread. In practice, modern CPU can easily handle 3 or 4 IO intensive threads per core, which leads to $L=12$ if one core is reserved for operating system only.

How many transactions does a MFN process and then has to store? Assume that $n$ transactions per second are uniformly divided among the $2^L-1$ concurrent blockchains. Then a MFN has to validate and store 
\begin{equation}
s_{MFN}=\frac{nL}{2^L-1}
\label{eq:num_of_trans_to_store_mfn}
\end{equation}
transactions per second and needs to download for each of them at most one chain composed of maximally $L$ sub-blocks for validation purposes, see Section \ref{sec:sharding_validation_and_storage}. A minimal amount of communication is achieved if each of those sub-blocks contains just one transaction. For this optimal case we obtain that the number of transactions per second one MFN has to download is bounded by
\begin{equation}
d_{MFN} = \frac{nL^2}{2^L-1} + L.
\label{eq:num_of_trans_to_download_mfn}
\end{equation}

\begin{figure}[ht]
	\scalebox{0.85}{\input{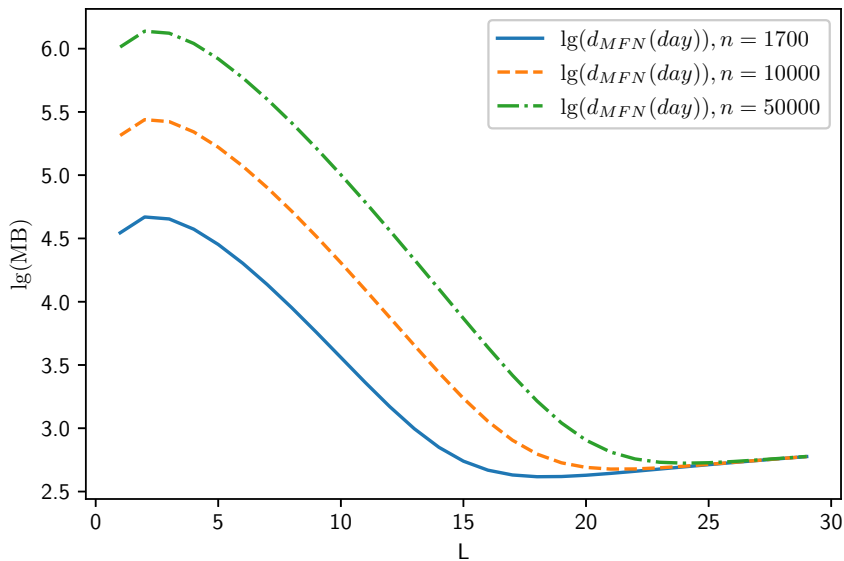}}
	\caption{An upper bound of minimal download per day in MB of a single MFN: result of running Algorithm \ref{alg:download_mfn}}
	\label{fig:download_mfn}
\end{figure}

Clearly, the requirement for storage \refe{eq:num_of_trans_to_store_mfn} are monotonically decreasing function of $L\ge 1.$ In Figure \ref{fig:download_mfn}, there is estimated daily download in megabytes for different values of $n$, starting from $1700$t/s which is close to the average number of VISA network transactions per second. This figure is generated by running Algorithm \ref{alg:download_mfn} which of course employs relation \refe{eq:num_of_trans_to_download_mfn}. We assume an average transaction size of $250$B. We see that to minimize communication, a certain minimal $L$ has to be reached, such that the linear term in \refe{eq:num_of_trans_to_download_mfn} starts to dominate. 

\begin{figure}[ht]
	\scalebox{0.85}{\input{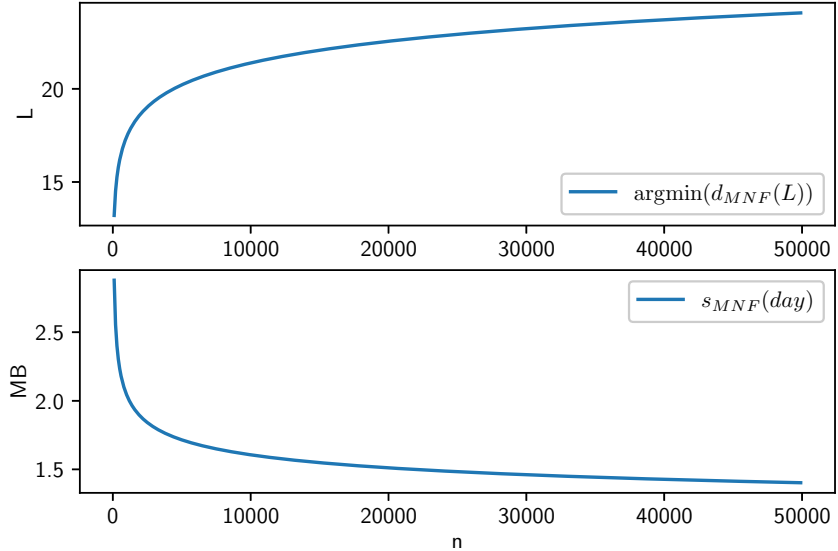}}
	\caption{$L = \mathrm{argmin}(d_{MFN}(L, n))$ for different transaction throughputs per second $n$ and corresponding storage per one MFN per day: result of running Algorithm \ref{alg:optimal_L_mfn}}
	\label{fig:optimal_L_mfn}
\end{figure}

Moreover, for each $n\in\NN^+$ there exist one minimum of \refe{eq:num_of_trans_to_download_mfn}. By running Algorithm \ref{alg:optimal_L_mfn} we compute those optimal $L$ for different transaction throughputs $n$. And we compute the corresponding storage per day for single MFN using \refe{eq:num_of_trans_to_store_mfn}. The result is presented in Figure \ref{fig:optimal_L_mfn}. By employing those $L$, we obtain per day downloads in range from approximately $300$MB for $L = 13$ to approximately $530$MB for $L=24$. It means that theoretically a 52kbits/s connection is sufficient (neglecting upload). Moreover, the per day MFN storage is monotonically decreasing function of $n$ for $L$ values minimizing the connection throughput, dropping from merely $\sim 2.9$MB for $L=13$ to $\sim 1.4$MB for $L=24$. 

The scaling potential should be now obvious. It could be possible to run a MFN node for $L=13$ which needs a little bit more then 1GB of transaction storage per year and can run on a 2G wireless connection. By scaling further, the storage actually drops to approximately $500$MB per year for $L=24$ and the 2G connection is still sufficient, neglecting latency issues.

\section{Conclusion}

In this paper, we have argued that a solution to the problem of the Bitcoin on-chain scaling lies outside of the comfort zone of the current concept of its blockchain, see e.g. Appendix \ref{app:mining_economics}. Nowadays, transactions requiring completely different security guaranties, considering the monetary values they move, are included into the same block. They are signed by the aggregate hashrate of all the participating miners and validated and stored by all the full nodes. Within these constraints, a scaled up Bitcoin blockchain simply cannot be found.

We have proposed a multi-scale approach, where extended Bitcoin transactions $t$ having two tunable parameters are uniquely assigned to a certain shard $\block_{l,s}, s\in\{0,\dots, 2^l-1\}$ of the network at a certain level $l\in\isetl.$ The default level $l$ assigned to a transaction is controlled by its output value per bit $\vpb$, but $l$ can be freely chosen by the user. The lower $l$, the higher the number of hashes used to secure $t$ and more expensive it becomes to double spend. On the contrary, both electric energy consumption and fees decrease as $l$ grows.

The sub-blocks $\block_{l,s}$ are organized in a (binary) tree topology with coordinates $(l,s)$. The transactions $t$ are uniquely assigned to individual nodes/branch\-es of this tree by a shard function $\shardf$, see \refe{eq:shard_function_mining_transactions}. This function is an essential concept, which facilitates probabilistic assignment of not only transactions $t$ to nodes $(l,s)$ of the tree, but also miners and other participants to those nodes, see \refe{eq:shard_function} and Section \ref{sec:sharding_validation_and_storage}. Metaphorically, transactions meet their miners and minimal full nodes (MFNs) at coordinates $(l,s)$ with uniform probabilities in $s$. A MFN is an updated version of full node. It validates and stores only transactions along its branch $(s_0, s_1, \dots, s_{L-1})$ assigned to it by \refe{eq:shard_function}. Consequently, the state $s$ is distributed among all MFNs. 

We have suggested that an efficient information exchange scheme among the participants should fully utilize the (binary) tree topology. In Section \ref{sec:hbs_binary_tree} this topology is strictly enforced, see Figure \ref{fig:hbs_tree}, but in Section \ref{sec:hbs_concurrent_levels} only loosely, see Figure \ref{fig:hbs_tree_concurrent}. In the latter, we essentially have $2^L+1$ interlocked concurrent blockchains. This constitutes our final scaling proposal which achieves all objectives. It offers thousands of timely  transactions per second with a variable level of trust together with aligned network fees and energy usage. And the state $s$ is distributed in a way which allows practically anybody to participate and run a MFN.

Many results in this paper are based on a simple security model \refe{eq:def_security} together with Assumption \ref{ass:same_security}. While \ref{eq:def_security} is certainly only a rough approximation of the reality and the expectations of some or even all the users can differ from those in Assumption \ref{ass:same_security}, it is of no concern. The solution found in Section \ref{sec:hbs_concurrent_levels} is sufficiently flexible and self correcting. First, the users are free to choose an appropriate security level $l\in\isetl.$ Second, miners may or may not include those transactions into their proposition blocks at the fees offered by the users. A true equilibrium between costs and guaranties will be found by market forces.

The Bitcoin elites, at first sight, have no economic incentive to support the changes proposed in this publication that include aiming for fees proportional to transaction value $v(t).$  Some of them are certainly satisfied with the current role of BTC as a store of value, volatile but a store of value nevertheless. But other crypto plutocrats may recognize that without resolving inconsistency between the costs and fees as block rewards decrease, not only will $\bitcoinA$ never become a universal mean of payment, but it will become progressively insecure and eventually at risk of becoming extinct, see \cite{Hasu2019}. The implications for the store of value proposition are obvious. 

Actually, mathematically speaking, there exists at least one very early adopt\-er that could find this contribution significant. The one who wrote:
\emph{\\``Craig Steven Wright is a liar and a fraud. He doesn't have 
the keys used to sign this message.\\
\\
The Lightning Network is a significant achievement. However, 
we need to continue work on improving on-chain capacity.\\
\\
Unfortunately, the solution is not to just change a constant 
in the code or to allow powerful participants to force out others.\\
\\
We are all Satoshi''}
  
We are highly optimistic about the Bitcoin future. The community will certainly recognize the economic discrepancies described in this paper and build upon the solutions proposed. In the end, Bitcoin ecosystem is as democratic as one can hope it to be. If developers suggest solutions and users choose to support them, the elites have to bend the knee or risk their fortune to evaporate, see Section 2.4 of \cite{Hasu2019}. There is more then enough time for discussions, solutions and implementations of them. But complacency and postponement are never good tendencies. It is better to start right now or others will take the lead.

\section{Acknowledgement}

I have never participated in Bitcoin's development or in any other crypto based project. This paper is written in my limited spare time and for prospective public good. If you liked what you read, you may donate and support my financial freedom to think and create. Please, do not be ashamed to donate any small contribution. In the end, $10^6\cdot100 sats = 1 \bitcoinA.$ Thank you.

Before sending any coins, please double check that the addresses here are the same as the ones published on my github account \cite{github} or check the sha256sum of this paper there. To communicate with a limited trust, please use my email\footnote{``Total paranoia is just total awareness.'' -- Charles Manson}.
\\

\noindent SegWit: bc1qrfunay7uekg5aj578suh6c6n8d9lqhwwf5h3ps\\
Legacy: 1EZRjWCagpegLcydVPSB4zFPgEdZ4wYuKu

\section{Note on republishing of paper}
\label{sec:republishing}
I had multiple reasons for publishing the paper under a pseudonym. Firstly, I wanted the work to be evaluated on its own merits, separate from my scientific credentials attached to it. I believe that the identity of the mind behind a work shouldn't influence the judgment of its worth.

Secondly, the solutions presented in the paper could lead to bitcoin achieving the fully decentralized payment network status at scale as originally envisioned. This is a delicate proposition, clearly against various special interests trying to influence the bitcoin evolution in their favor. For instance, if the base layer of bitcoin remains inefficient, there's a necessity for separate layer 2 solutions, potentially yielding substantial profits for those involved.

Lastly, publishing anonymously allowed me to relax certain aspects, such as conducting a less exhaustive literature review. It indeed appears evident that blockchain scaling has been analyzed in scientific literature from various perspectives, even beyond the realm of bitcoin. However, I did not have to concern myself with this and could naturally delve into the problem, addressing it comprehensively without fretting over details and losing focus. Ultimately, when starting from scratch with only an idea, it is practically impossible to obtain only previously published results. Usually, on the contrary, it can be a rather inventive path to take. Nevertheless, even though the paper does not aim for publication in a scientific journal, it is written rather carefully - ``This is the way.''

\begin{comment}Lastly, it was a social experiment. Can a single human mind be found within a reasonable final time interval that will appreciate the content so much that an anonymous author will be rewarded with few satoshis? It definitely failed. 

 So even if the certain quality was recognized, I was never awarded any finantial incentive that could potentially encourage me to build upon what I have done. Strikingly, in my view this paper represents potentially a very significant contribution to a network that has the current market capitalization of almost 1.3T\$.
\end{comment}

After two years, the paper was downloaded from viXra only a limited number of times. That's why it's even more surprising that I have received at least two proposals to publish the article in a scientific journal, one even being interesting. I do not know if the number is currently higher, since I have recently lost access to my anonymous email due to the inactivity window. This was the last straw. It might be simply by design that publishing a scientific paper anonymously is difficult. To sum up, I have decided to attach my name to this contribution because it allows me to disseminate it more easily and obtain some feedback. But, it was an interesting experiment.

Last but not least, after publishing my work, I have, more or less by coincidence, discovered that Peter Todd considered hierarchical block structure at least as far back as February 2014 \cite{Todd2014}. However, I could not find almost any published resources that provided sufficient details. Then, I discovered a podcast where he describes the idea in a little bit more depth \cite{Todd2014b}. Even if the presentation is informal and lacking rigor, the proposed structures are roughly the same. I was glad that I found those references since they supported my approach. Equally, I am glad that I found them after publication since I would probably never have done this analysis, knowing that the core idea was already proposed and moreover did not gain any significant traction. Why not? I have ``no'' idea.

\bibliographystyle{plain}
\bibliography{b2}

\begin{thebibliography}{10}

\bibitem{Back2002}
Adam Back.
\newblock {Hashcash - A Denial of Service Counter-Measure}.
\newblock http://www.hashcash.org/papers/hascash.pdf, 2002.

\bibitem{Beres2020}
Ferenc B\'eres, Istv\'an~A. Seres, and Andr\'as~A. Bencz\'ur.
\newblock {A Cryptoeconomic Traffic Analysis of Bitcoin's Lightning Network}.
\newblock {\em Cryptoeconomic Systems}, 2020.

\bibitem{Buterin2016}
Vitalik Buterin.
\newblock {A Proof of Stake Design Philosophy}.
\newblock
  https://medium.com/@VitalikButerin/a-proof-of-stake-design-philosophy-506585978d51.

\bibitem{Buterin2013}
Vitalik Buterin.
\newblock {Ethereum Whitepaper}.
\newblock {https://ethereum.org/en/whitepaper/}.

\bibitem{Carstens2021}
Agust\'{i}n Carstens, Jerome~H. Powell, Gillian Tett, and Jens Weidmann.
\newblock {How can central banks innovate in the digital age?}
\newblock https://www.bis.org/events/bis\_innovation\_summit\_2021/agenda.html,
  2021.

\bibitem{Clauson2021}
Aaron Clauson, Andrew Chow, Anthony Towns, Bruno Garcia, Fabian Jahr, fanquake,
  Hennadii Stephanov, Jon Atack, Luke Dashjr, MarcoFalke, Pieter Wuille,
  practicalswift, radymmcmillan, Sjors Provoost, Vasil Dimov, and W.J. van~der
  Laan.
\newblock {0.21.1 Release Notes}.
\newblock https://bitcoin.org/en/releases/0.21.1.

\bibitem{difficulty}
Bitcoin community.
\newblock Difficulty.
\newblock https://en.bitcoin.it/wiki/Difficulty.

\bibitem{BCH}
Wikipedia contributors.
\newblock {Bitcoin Cash --- Wikipedia, The Free Encyclopedia}.
\newblock {https://en.wikipedia.org/wiki/Bitcoin{\textunderscore}Cash}, 2021.

\bibitem{Crain2016}
Brian~Fabian Crain and S\'ebastien Couture.
\newblock {EB65 -- Adam Back \& Greg Maxwell: Sidechains Unchained}.
\newblock
  {https://www.youtube.com/watch?v=jE{\textunderscore}telgnlw3M\&t=4775s}.

\bibitem{Decker2013}
Christian Decker and Roger Wattenhofer.
\newblock {Information Propagation in the Bitcoin Network}.
\newblock In {\em {13-th IEEE International Conference on Peer-to-Peer
  Computing}}, 2013.

\bibitem{Nah2004}
Fiona Fui-Hoom~Nah.
\newblock {A study on tolerable waiting time: how long are Web users willing to
  wait?}
\newblock {\em {Behaviour \& Information Technology}}, 23(3):153--163, 2004.

\bibitem{Grunspan2018}
Cyril Grunspan and Ricardo P\'erez-Marco.
\newblock {Double Spend Races}.
\newblock {\em {International Journal of Theoretical and Applied Finance}},
  21(08):1850053, 2018.

\bibitem{Hasu2019}
Hasu, James Prestwich, and Brandon Curtis.
\newblock A model for bitcoin's security and the declining block subsidy.
\newblock
  https://uncommoncore.co/wp-content/uploads/2019/10/A-model-for-Bitcoins-security-and-the-declining-block-subsidy-v1.02.pdf,
  2019.

\bibitem{Kerber2019}
Thomas Kerber, Aggelos Kiayias, Markulf Kohlweiss, and Vassilis Zikas.
\newblock {Ouroboros Crypsinous: Privacy-Preserving Proof-of-Stake}.
\newblock In {\em 2019 IEEE Symposium on Security and Privacy (SP)}, pages
  157--174, 2019.

\bibitem{BIP141}
Eric Lombrozo, Johnson Lau, and Peter Wuille.
\newblock {Segregated Witness (Consensus layer)}.
\newblock https://github.com/bitcoin/bips/blob/master/bip-0141.mediawiki.

\bibitem{Nakamoto2008}
Satoshi Nakamoto.
\newblock {Bitcoin: A Peer-to-Peer Electronic Cash System}.
\newblock https://bitcoin.org/bitcoin.pdf, 2008.

\bibitem{nullc2021}
nullc.
\newblock https://old.reddit.com/r/Bitcoin/comments/mtugta/mentor\_mon\-
  day\_april\_19\_2021\_ask\_all\_your\_bitcoin/gv86j6b/.

\bibitem{Rosenfeld2014}
Meni Rosenfeld.
\newblock {Analysis of Hashrate-Based Double Spending}.
\newblock https://arxiv.org/abs/1402.2009, 2014.

\bibitem{github}
Qui Somnium.
\newblock Hbs.
\newblock https://github.com/quisom/hbs, 2021.

\bibitem{Somnium2022}
Qui Somnium.
\newblock {Saving Proof-of-Work by Hierarchical Block Structure: Bitcoin 2.0?}
\newblock https://vixra.org/abs/2204.0040, 2022.

\bibitem{Todd2014b}
Peter Todd.
\newblock {LTB E104 - Tree Chains with Peter Todd}.
\newblock https://soundcloud.com/mindtomatter/ltb-e104-tree-chains-with.

\bibitem{Todd2014}
Peter Todd.
\newblock {Tree Chains}.
\newblock https://github.com/petertodd/tree-chains-paper.

\end{thebibliography}

\appendix 
\section{On Bitcoin, its economics, energy consumption and scaling}
\label{app:bitcoin_economics}

\subsection{In short, how Bitcoin works?}
Bitcoin \emph{transactions} are stored in \emph{blocks} of a predefined maximal size. Each subsequent block contains a hash of a previous block, thus proving that its creator accepts all the transactions recorded in this and all the earlier blocks. Consequently, the blocks form a \emph{blockchain}. \emph{Mining nodes} in the network compete to create a next valid block of transactions by finding in parallel and with a minimal amount of coordination a nonce for their proposition block such that its hash is smaller than a predefined number which determines the difficulty of finding such a hash. This difficulty is regularly adjusted, reflecting any change of the compound computational power of all the nodes, such that the block creation has a specified expected speed (e.g. one $1$ MB block per 10 minutes). When a hash for a certain block is found that fulfills the difficulty criterion and thus proving that the set amount of computational time has been spent to secure the transactions contained in the block, this block is added to the blockchain. The creator of the block then obtains a certain amount of newly minted bitcoin called \emph{block reward}, this being the main incentive to take part in signing of the transactions.

\subsection{Deflationary nature}
Some argue that Bitcoin's technological constraints such as its scalability issues are a secondary problem to its ultimately deflationary nature. They claim that bitcoin's appreciation in market value discourages people from using it for payments. It is a reasonable hypothesis, since, its finite supply is caped at 21M and ultimately it should be a deflationary coin. Its inflation rate is governed by block reward which halves every 210000 mined blocks, i.e. about every 4 years. It was initially 50 bitcoins (BTC), currently being 6.25 BTC. Bitcoin was programmed to be better in scarcity than gold. And there are other aspects that make it superior to gold, e.g. its divisibility or its zero mass which implies practically zero costs holding the asset. Even the most influential people in finance now recognize or at least admit that BTC shows some gold-like properties and behavior. Let us quote Federal Reserve Chairman Jerome Powell talking about cryptoassets and bitcoin in particular: \emph{``\dots There are also not particularly in use as means of payment. It's more a speculative asset that's essentially a substitute for gold, rather than for the dollar.\dots''} \cite{Carstens2021}. 

Inflationary fiat/crypto currencies are definitely easier to be disposed of when paying for goods or services. If BTC's value rises when it is still inflationary due to the ongoing minting process, what will happen when this stops? 

%Yes, it seems obvious that its halving scheme is pyramidal by construction. Arguably it is a good thing since the dropping inflation rate is crucial for Bitcoin's mass adaption. 

But rather than focusing on the question whether BTC is the new gold, we focus in this paper on what Bitcoin was meant to be and it clearly has not yet become: a successful general payment system. Without removing the existing technological hurdles, we will never really know if people are or will ever be prepared to use BTC for daily payments and in this way to potentially spread real wealth to others, maybe being even ``slightly'' altruistic. In the end, it is difficult to assess whether a deflationary monetary system could work, especially because societies are used to practically constant presence of inflation. Philosophically speaking, human life time is similarly deflationary. Each day we have one day less of it.

\subsection{Current state of blockchain affairs}
\label{sec:current_state_blockchain_affairs}

There exist more than 5000 blockchain based cryptoassets at the time of writing this article, precise number not being important. Most of them are merely \emph{tokens} utilizing one of a handful of public blockchain implementations which support smart contracts. The pioneering one is Ethereum \cite{Buterin2013}, the number two player with still less than a half of Bitcoin's market cap. They together represent more than $60\%$ of more than $2$ trillion USD invested in public blockchains.

The Ethereum project extends in a very meaningful way Bitcoin's codebase and builds a truly programmable blockchain. Bitcoin itself is a virtual machine evaluating small programs but limited by design. Arguably, its Taproot upgrade from November 2021 expands substantially its smart contract abilities, potentially reducing Ethereum's significant usefulness advantage in this domain \cite{Clauson2021}. Both Bitcoin and Ethereum use PoW based consensus algorithms but Ethereum ecosystem is long planned to switch to \emph{Proof of Stake} (PoS) based one \cite{Buterin2016}.

Proof of Stake based consensus algorithms, among others, are proposed as alternatives to PoW to mitigate some of its shortcomings, mainly its negative ecological impact. In short, the blocks are created by validators that are chosen randomly with probability proportional to their staked holdings of the underlining coin. They are rewarded for creating a valid block and/or punished for misbehaving by loosing a part of or the whole stake.  The idea is simple and intriguing at first but a complex one to apply at a closer look: see e.g. \cite{Crain2016} for critical counterarguments. It suffices to note that the Ethereum network is still running on PoW even after more than 5 years of PoS development.

Among PoS based blockchain platforms, Cardano clearly stands out. Its market capitalization rose from roughly 1 billion USD in 2020 to approximately 100 billion USD in 2021, becoming the number tree player in the public blockchain market. This may be to a great extent a speculation move into a limited smart contract enabled PoS blockchain space, driven by the rising need for an environmentally friendly blockchain technology. Nevertheless, Cardano has in our view a unique selling proposition since its development is guided by evidence-based methods and its consensus algorithm Ouroboros is formally verified under different models  in a sequence of peer-reviewed scientific papers, see \cite{Kerber2019} and the references wherein.\footnote{It is up to a competent reader to independently evaluate whether the models and their assumptions are realistic.}

The practical evidence that PoS can be implemented robustly enough to secure values is indeed comparatively limited. PoW is much more battle tested and clearly wins with respect to its simplicity which helps to build trust. Moreover, Bitcoin and Ethereum are much more then just their consensus algorithms, supporting in return the PoW momentum. Only time will tell whether PoS wins similar levels of market approval, this almost certainly being dependent on the success of Ethereum's transition to PoS. 

%ADA, Cardano's coin was distributed via ICO and the major beneficiaries have currently no %reasons to kill its meteoric rise. 
%Sybil resistance of Bitcoin's PoW.

% accumulated trust and values it secures for years, and most importantly, its theoretically provable guaranties. 

%But incentives to present a viable PoS alternative are significant.

\subsection{Bitcoin energy consumption}
\label{app:energy}

The power efficiency $\pe$ of modern commercially available ASIC miners is approximately $30$ J/TH. Given the current approximate total hashrate $\thr$ of $1.20e8$ TH/s, corresponding to $1.2$M of modern Antminers S19j Pro, we arrive at the power of  $3.6$ GW, equivalent to $3-4$ conventional nuclear rectors.

It is equivalent to $3.6e6$ kWh. At electricity price $\ep$ e.g. $0.1$\$ per kWh, we arrive at one hour energy costs of $360000$\$. Since one block is mined on average every 10 minutes, energy costs per block are cca $60000$\$. For a legacy transaction which size $\bsize(t)$ equals to $250$B, costs are

\[
\approx \frac{\displaystyle 250\mbox{B}}{\displaystyle 1024^2\mbox{B}} \cdot 60000 \mbox{\$} = 14.3\mbox{\$}
\]
if blocks are assumed to be fully filled. Altogether, a general formula for transaction electric energy consumption $\tec$ is 

\begin{equation}
\label{eq:tec}
\tec\mbox{[kWh]} := \pe[\mbox{J/TH}] \cdot\thr\mbox{[TH/s]}\cdot\frac{600\mbox{s}}{3.6e6}
\cdot\frac{\displaystyle \bsize(t)\mbox{[B]}}{\displaystyle 1024^2\mbox{B}}
\end{equation}
and the cost per transaction $cpt$ is then simply
\[
\cpt \mbox{[\$]} := \tec[\mbox{kWh}]\cdot\ep\mbox{[\$/kWh]}.
\]

Certainly, these are lower estimates, since not all the hardware used to mine bitcoin is recent and it does not include all the indirect energy use. Nevertheless, this simple estimate yields that Bitcoin network uses approximately $32$TWh of electricity per year - far less than normally estimated. 

But indeed, because of the block reward, even mining on older hardware is profitable. We may obtain an upper estimate of Bitcoin energy use, if we know the average price of electric energy used for mining and assume that miners are rational and stop mining when it starts to be unprofitable.

Let $\mbox{BTCUSD}$ be the bitcoin exchange rate with USD, $\tfpb$ the transaction fee per bit and $\br$ the block reward. The total block reward $\tr$ is

\[
\tr\mbox{[\$]} = \br[\bitcoinA] \cdot \mbox{BTCUSD} + 1024^2\mbox{B} \cdot 8\tfpb\mbox{[\$/b]}.
\]

For a rational miner 

\begin{equation}
\tr\mbox{[\$]} >  \ep\mbox{[\$/kWh]} \cdot \epb\mbox{[kWh/block]},
\label{eq:totalreward_energy_inequality}
\end{equation}
where $\epb$ is the total electrical energy consumption of the Bitcoin network to mine a block which we estimated above as

\[
\epb \approx  \pe[\mbox{J/TH}] \thr\mbox{[TH/s]}\frac{600\mbox{s}}{3.6e6}.
\]

From inequality \refe{eq:totalreward_energy_inequality} we get an upper bound on electric energy consumption based on the market prices. Currently the $\tfpb$ is stabilized under $0.001875$ \$/b and BTCUSD is above 40k. Consequently, the network total consumption should not be higher then
\[
\epb\mbox{[kWh/block]} < \frac{6.25\cdot 40000 + 8 \cdot 0.001875 \cdot 1024^2}{0.1} \approx 2.66 \mbox{GWh}
\]
for the electricity price of $0.1$ \$/kWh. This leads to the upper bound of approximately $140$ TWh per year. This is much closer to the other estimates. Nevertheless, we think that the real consumption is much lower, since electric energy is the main cost of mining and professional miners have all incentives to use as efficient hardware as possible.

It is important to note that transaction costs are normally relatively small with respect to the block reward. Even during its historical maximum of $\tfpb\approx 0.03$ \$/b, the block reward was a bigger chunk of total reward, even if comparable.

As the block reward (i.e. the inflation rate) of bitcoin will decrease, the electricity consumption should be governed more end more by transaction fees only and per transaction it should be at least a few magnitudes lower, if we want the Bitcoin network to succeed.

%6.25*40000 + 1024^2/250*3

\subsection{Mining economics and scaling}
\label{app:mining_economics}

%Let us now consider a different perspective, that of miners. 

Mining is a very competitive business and we assume that the majority of miners does not behave altruistically but they rationally maximize their profits. Their financial incentive to sign transaction is composed of two parts: block rewards and transaction fees. 
	
The cost structure of bitcoin mining is complex. Nevertheless, except all the initial investments to start a mining business, including the equipment acquisition costs, the electric energy consumption is the main production cost, e.g. denominated in USD.
	
To enlighten the peculiar nature of Bitcoins minting scheme we can employ a simple mind experiment: Let us imagine for a moment we fix the whole Bitcoin ecosystem including the energy consumption and its USD denominated price $\ep\mbox{[\$/kWh]}$, only halving the block reward in BTC. For miners to recuperate the block reward loss, the transaction fees and/or the market price of BTC have to increase.
	
Let us now further assume that the BTC market  price will stagnate at preceding levels after such a hypothetical halving. Then there exists an immediate pressure to increase the transaction fee per byte to recuperate the USD denominated loss.
	
Conversely, let us now assume that the market price of BTC rises proportionally after this hypothetical halving, i.e. it at least doubles since the inflation rate is halved. Then there is no increased pressure to hike the transaction fees, since the USD dominated reward of miners is at least preserved. Nevertheless, the transaction fee at least doubles in USD denominated terms.
	
In both cases, the halving is making bitcoin payment system more expensive and thus less competitive. In reality, a mix of the above two extremal possibilities happens. Both the BTC market price rises, even if very non-linearly, and the transaction fees rise as well.  Moreover, eventual more than doubling of BTC market price incentives old miners to increase capacity and new miners to join, increasing the overall energy consumption of the Bitcoin network and rising the environmental concerns.
	
If the transaction throughput rose at least proportionally, i.e. if it at least doubled after each halving, it could be argued that real costs per transaction are preserved. But it would be desirable that utilization increases quicker or even dramatically, such that the energy usage per transaction $\tec$, estimated in Equation \refe{eq:tec}, could drop significantly as well.
	
According to Equation \refe{eq:tec}, $\tec$ can be reduced by reducing the total hash rate $\thr$, by increasing the power efficiency $\pe$ of mining equipment, by increasing the block generation speed, by decreasing the transaction size or finally by increasing the block size. Let us shortly analyze the individual possibilities. 
	
First, any uneven, substantial reduction of $\thr$ theoretically negatively impacts Bitcoin's network security, since potential attack abilities of some miners increase. This is certainly not desirable, so if at all possible, only a smooth, even and slow $\thr$ decrease could be promoted. 
	
On the other hand, incentives to increase the power efficiency $\pe$ of mining hardware are innate and $\pe$ also rises at a steady pace. 
	
One successfully implemented \emph{Bitcoin Improvement Proposal} (BIP) which inter alia also decreased transaction size is Segregated Witness \cite{BIP141}. 
	
Further, the most direct way how to decrease $\tec$ and at the same time to allow for quick/instant payments is to generate new blocks more frequently. Moreover, it is proven practically possible by Ethereum's implementation of PoW, where new blocks are generated approximately once every 15 seconds \footnote{They are however much smaller.} This would only require to distribute the reward into multiple block inside the current 10 minute interval, so that the Bitcoin's inflation rate is preserved.
	
The last obvious possibility is to increase maximum size of blocks. It is a very sticky issue on which some very intelligent people disagree \cite{BCH}. In short, the supporters of a block size increase emphasize an on-chain medium of exchange function. The main argument against is that larger blocks would make full nodes more expensive to operate and this would lead to further centralization which weakens Bitcoin's values proposition. Actually, the same counterargument can be employed against the above mentioned increase in block generation frequency.
	
We agree that the decentralized nature of Bitcoin is the most sacred property of Satoshi Nakamoto's invention and should be preserved at all costs. We explain in the next Section that even if centralization of mining is partially inevitable, it is indeed crucial to promote running full Bitcoin nodes and thus any substantial increase in block size or generation speed is problematic. 

Inevitably, any viable scaling idea has to represent a substantial evolution of blockchain technology as we know it. One such an off-chain approach is Lighting Network \cite{Beres2020}. In this paper we propose an on-chain solution. 

\subsection{On centralization of mining}
\label{app:centralization_of_mining}

First, we postulate that centralization is practically inevitable in any system where an increase in efficiency leads to a corresponding increase in profits. This is certainly true for bitcoin mining. The actors who are first able to exploit the inefficiencies of such a system, acquire more wealth which in turn allows then to extract more of the inefficiency. As inefficiency potential is diminishing in time, the first players have an enormous advantage. 

In natural real world systems, the external conditions often change, so new inefficiencies and opportunities to exploit them emerge. Bitcoin is however a man-made system governed by strict rules of its protocol which can be changed only by consensus and thus rather stable and consequently it supports inertia, resulting naturally in centralization of mining, mainly due to economies of scale.

This kind of centralization is not a problem in isolation, even if it seems to be the case according to classical Byzantine fault tolerance (BFT) analysis, established for Bitcoin's consensus protocol already in \cite{Nakamoto2008}. Satoshi Nakamoto however only estimated the probability of successful double-spending attack. For a more detailed analysis see \cite{Rosenfeld2014}, where the correct probabilistic model is employed. Moreover, the paper briefly discusses economics of double spending as well. Finally, a rigorous mathematical treatment of Bitcoin's BFT probabilistic model is presented in \cite{Grunspan2018}. 

However, the bitcoin mining is already centralized to such an extend that BFT is insufficient to ensure its security. E.g. only the five biggest known mining pools control together more than $50\%$ of hashrate. BFT implies that a majority attacker can easily double-spend, censor transactions and even claim all block rewards for himself. Had it been rational to execute these attacks, it would have been rather easy for these players to organize. But it is against their best interest as already recognized in \cite{Nakamoto2008}:\emph{``If a greedy attacker is able to assemble more CPU power than all the honest nodes, he would have to choose between using it to defraud people by stealing back his payments, or using it to generate new coins. He ought to find it more profitable to play by the rules, such rules that favour him with more new coins than everyone else combined, than to undermine the system and the validity of his own wealth.''} 

This observation is well established in \cite{Hasu2019}, where a simple but therefore a robust cost-reward based Bitcoin security model is presented. The main conclusion is that miners are committed to preserve trust in Bitcoin network, since their long-term infrastructure investments are non-repurposable due to their high level of specialization. Any loss of trust normally results in BTC price plummeting, affecting the margins significantly. Following the estimates in the paper, at the current issuance of 6.25 BTC per block only a 5$\%$ sustainable market price drop would result in almost a 10K BTC loss for an attacker with $60\%$ hashrate majority.

We point out that the authors in \cite{Hasu2019} did not analyze a sophisticated attacker which hedges against BTC market drop by e.g. a short position in bitcoin futures. However, such an attack is prohibitively expensive since e.g. Interactive Brokers apply for short bitcoin futures (BRR) positions margins of $150\%$  of daily settlement price. Moreover the liquidity of bitcoin futures is limited with open interest currently around 1.5 billion USD. If deep bitcoin derivatives markets develop, such attacks become more probable. Arguably, Bitcoin's ecosystem inclusion into the current financial system may represent its biggest existential threat.

So if the mining centralization is currently not an existential problem, which type of decentralization is necessary to be preserved? Blocks created by the miners are ultimately validated by Bitcoin users running full nodes. If a majority of these users is against the blocks proposed by the miners or their behavior in general, they can organize off-chain to suspend Nakamoto consensus. For details and precedents from Bitcoin's history, see again  \cite{Hasu2019}. This is the ultimate circuit breaker. Imagine e.g. a last-resort capital punishment scenario when the users would change the hashing algorithm ``overnight'', rendering the special mining equipment practically immediately useless. This ability has to be preserved at all costs. It requires that many users have to be able to run full-node Bitcoin clients. More precisely, the Bitcoin users have to be able to validate transactions included into the blocks. This together with the ability of users to sell their bitcoins, secures a healthy functioning of Bitcoin's network.

\section{Algorithms}
\label{app:algorithms}

\begin{algorithm}
\caption{Hashing time scales linearly with size.}
\label{alg:hashing_time}
\begin{verbatim}
import hashlib
import time
import os
import matplotlib.pyplot as plt
plt.rcParams['pgf.texsystem']= 'pdflatex'

n = 100
x = range(n)
y = []
for i in x:
    m = os.urandom(i*1024**2)
    time_b = time.time()
    hashlib.sha256(m).hexdigest()
    time_e = time.time()
    y.append(time_e - time_b)

fig = plt.subplots()    
plt.plot(x,y)
plt.xlabel("size[MB]")
plt.ylabel("time[s]")
plt.savefig('hash_time.pgf')
\end{verbatim}
\end{algorithm}

\begin{figure}[h]
	\scalebox{0.85}{\input{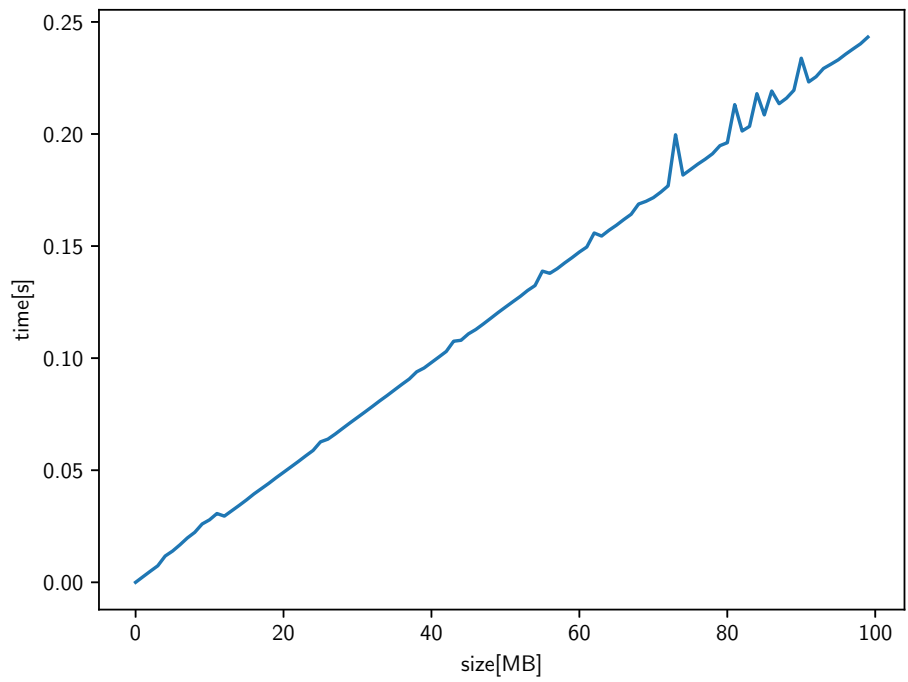}}
	\caption{Linearity of hashing time: result of running Algorithm \ref{alg:hashing_time}}
	\label{fig:hash_time}
\end{figure}

\begin{algorithm}
\caption{Prepares Bitcoin transaction dataset for analysis.}
\label{alg:tx_dataset}
\small
\begin{verbatim}
import os
from blockchain_parser.blockchain import Blockchain
import random
import pandas as pd

columns_blk = ['block_height', 'block_hash', 'version',\ 
    'previous_block_hash', 'merkle_root', 'timestamp',\ 
    'bits', 'nonce', 'difficulty']
types_blk = ['uint32', 'str', 'uint32', 'str', 'str',\
    'uint32', 'uint32', 'uint32', 'float']
dtypes_blk = {c: t for (c, t) in zip(columns_blk, types_blk)}
columns_tx = ['block_height', 'txid', 'hash', 'version', 'n_inputs',\
    'n_outputs', 'is_segwit', 'is_coinbase', 'size', 'output_value']
types_tx = ['uint32', 'str', 'str', 'uint32', 'uint32',\
    'uint32', 'bool', 'bool', 'uint32', 'uint64']
dtypes_tx = {c: t for (c, t) in zip(columns_tx, types_tx)}
blk_df = pd.DataFrame(columns = columns_blk)
blk_df = blk_df.astype(dtype = dtypes_blk)
tx_df = pd.DataFrame(columns = columns_tx)
tx_df = tx_df.astype(dtype = dtypes_tx)

blockchain = Blockchain(os.path.expanduser('.../btc/blocks'))
bx_list = []
for blk in blockchain.get_ordered_blocks(\
    os.path.expanduser('.../btc/blocks/index'),\
    start = 650000, end = 700000):
    if (random.random() > 0.01):
        continue
    h = blk.header
    df = pd.DataFrame([[blk.height, blk.hash, h.version,\
    h.previous_block_hash, h.merkle_root, h.timestamp,\
    h.bits, h.nonce, h.difficulty]], columns = columns_blk)
    blk_df = blk_df.append(df)
    tx_list = []
    for tx in blk.transactions:
        sum_out = 0
        for tout in tx.outputs:
            sum_out = sum_out + tout.value
        tx_list.append([blk.height, tx.txid, tx.hash, tx.version,\
        tx.n_inputs, tx.n_outputs, tx.is_segwit, tx.is_coinbase(),\
        tx.size, sum_out])
tx_df = tx_df.append(pd.DataFrame(tx_list, columns = columns_tx))

tx_df.reset_index(drop = True, inplace = True)
blk_df.reset_index(drop = True, inplace = True)
tx_df.to_hdf('transactions.h5', key='tx_df', mode='w')
blk_df.to_hdf('transactions.h5', key='blk_df')
\end{verbatim}
\end{algorithm}

\begin{algorithm}
\caption{Algorithm to plot graph in Figure \ref{fig:tx_vpb}.}
\label{alg:plot_tx_hist}
\begin{verbatim}
import numpy as np
import pandas as pd
import matplotlib.pyplot as plt
plt.rcParams['pgf.texsystem']= 'pdflatex'

T = pd.read_hdf('transactions.h5', 'tx_df')
T = T[T['output_value'] > 0]
vpb = T['output_value']/(8*T['size'])

count, bins, ignored = plt.hist(np.log10(vpb),\
bins = 100, density = True)

mu = np.sum(np.log10(vpb))/len(vpb)
sigma = np.sqrt(np.sum((np.log10(vpb) - mu)**2)/(len(vpb)-1))

x = np.linspace(min(bins), max(bins), 100)
pdf = (np.exp(-(x - mu)**2 / (2 * sigma**2))
/ (sigma * np.sqrt(2 * np.pi)))
plt.plot(x, pdf, linewidth = 1.5, color = 'r')
plt.xlabel(r'\lg($\beta(t)$)')
plt.ylabel('density')
plt.savefig('tx_vpb.pgf')
plt.show()
\end{verbatim}
\end{algorithm}

\begin{algorithm}
\caption{A Python implementation of Algorithm \ref{alg:tx_segmentation}.}
\label{alg:tx_segmentation_python}
\begin{verbatim}
import numpy as np
import pandas as pd
	
def segment_transactions(L, T, key):
    T = T.sort_values(ascending = False, by = key)
	
    M = T[key].values[0]
    m = T[key].values[-1]
	
    S_L = (np.log10(M)-np.log10(m))/L
    #S_L = (np.ceil(np.log10(M))-np.floor(np.log10(m)))/L
	
    indx = [0]*(L+1)
    indx[0] = 0
    s = T[key]
    C = np.log10(M) - S_L; l = 0
    for r, v in enumerate(s):
        if np.log10(v) < C:
            l = l+1
            indx[l] = r
            C = C - S_L
    indx[L] = len(s)
    Ts = []
    for l in range(L):
        Ts.append(T[indx[l]:indx[l+1]])
    return Ts
if __name__ == "__main__":
    tx_df = pd.read_hdf('transactions.h5', 'tx_df')
    tx_df.reset_index(drop = True, inplace = True)
    tx_df = tx_df[tx_df['output_value'] > 0]
    tx_df['vpb'] = tx_df['output_value']/(8*tx_df['size'])
    Ts = segment_transactions(6, tx_df, 'vpb')
\end{verbatim}
\end{algorithm}

\begin{algorithm}
	\caption{Algorithm to obtain latex code of Table \ref{table:segmented_txs}.}
	\label{alg:get_table_segmented_txs}
	\small
	\begin{verbatim}
def generate_latex_table(Ts):
    text = """\def\\arraystretch{1.1}
\\begin{table}
\\begin{tabular}{|"""
    text += 'c|'*(len(Ts)+1)+"""}
\hline
$l$"""
    for i in range(len(Ts)):
        text += f' & {i}'
    text +="\\\ \hhline{|"
    text += '=|'*(len(Ts)+1)+"""}
$|\mathcal{T}_l|$""" 
    for i in range(len(Ts)):
        text += f' & {len(Ts[i])}'
    text += """ \\\ \hline
min($\\beta(t)$)""" 
    for i in range(len(Ts)):
        text += f' & {Ts[i]["vpb"].min():.2}'
    text += """ \\\ \hline
max($\\beta(t)$)""" 
    for i in range(len(Ts)):
        text += f' & {Ts[i]["vpb"].max():.2}'
    text += """ \\\ \hline
$\overline{\\beta(t)}$""" 
    for i in range(len(Ts)):
        text += f' & {Ts[i]["vpb"].mean():.2}'
    text += """ \\\ \hline
min($v(t)$)""" 
    for i in range(len(Ts)):
        text += f' & {Ts[i]["output_value"].min():.2}'
    text += """ \\\ \hline
max($v(t)$)""" 
    for i in range(len(Ts)):
        text += f' & {Ts[i]["output_value"].max():.2}'
    text += """ \\\ \hline
$\overline{v(t)}$""" 
    for i in range(len(Ts)):
        text += f' & {Ts[i]["output_value"].mean():.2}'
    text += """ \\\ \hline
$\sum v(t)$""" 
    for i in range(len(Ts)):
        text += f' & {Ts[i]["output_value"].sum():.2}'
    text += """ \\\ \hline
\end{tabular}
\caption{...}
\label{...}
\end{table}"""
    return text
	\end{verbatim}
\end{algorithm}

\begin{algorithm}
	\caption{Functions to estimate $\cti$, $\{\ti_l\}$ and time investment per level for our example dataset. Input Ts is obtained by Algorithm \ref{alg:tx_segmentation_python}.}
    \label{alg:compute_c_eta_and_eta}
	\begin{verbatim}
def compute_c_eta_and_eta(Ts, num_blocks = 493):
    s = 0
    for l in range(len(Ts)):
        s += Ts[l]['vpb'].mean()*8*Ts[l]['size'].sum()
    c_eta = 600*num_blocks*10**8/s
    eta = []
    for l in range(len(Ts)):
        eta.append(c_eta*Ts[l]['vpb'].mean()/10**8) 
    return c_eta, eta

def compute_time_per_level(eta, Ts, num_blocks = 493):
    t = []
    for l in range(len(Ts)):
        t.append(eta[l]*8*Ts[l]['size'].sum()/num_blocks)
    return t
	\end{verbatim}
\end{algorithm}

\begin{algorithm}
\caption{A python implementation of shard function \refe{eq:shard_function}.}
\label{alg:shard_func}
\begin{verbatim}
def shardf(l, string):
    b = hashlib.sha256(bytes(string, 'utf-8')).digest()
    shifts = []
    for r in range(l//8+1):
        shifts = shifts + [int(i) for i in f'{b[r]:08b}']
    shards = [0]*(l+1)
    for i in range(l):
        if shifts[i] == 0:
            shards[i+1] = 2*shards[i]
        else:
            shards[i+1] = 2*shards[i]+1
	return shards[l], shards
\end{verbatim}
\end{algorithm}

\begin{algorithm}
\caption{Algorithm to plot graph in Figure \ref{fig:sharding_efficiency}.}
\label{alg:sharding_efficiency}
\begin{verbatim}
import numpy as np
import matplotlib.pyplot as plt
plt.rcParams['pgf.texsystem']= 'pdflatex'

def ratio_blocks_mfn_vs_full(N, L):
    return N*L**2/(2**L-1)**2 + L/(2**L-1)

n_trans = 4200
L = 25
x = range(1, L)

y1 = [np.log10(ratio_blocks_mfn_vs_full(n_trans,i)) for i in x]
y2 = [np.log10(ratio_blocks_mfn_vs_full(2**i-1,i)) for i in x]
y3 = [np.log10((2**i-1)/600) for i in x]

fig, ax = plt.subplots()    
ax.axhline(y = 0, color = 'k', linewidth = 1)
ax.plot(x, y1, label = '$\lg(r), N = 4200$', linestyle = '-')
ax.plot(x, y2, label = '$\lg(r), N = 2^L-1$', linestyle = '--')
ax.plot(x, y3, label = '$\lg((2^L-1)/600)$', linestyle = '-.')
ax.legend()
plt.xlabel('L')

plt.savefig('sharding_efficiency.pgf')
\end{verbatim}
\end{algorithm}

\begin{algorithm}
\caption{Algorithm to plot graph in Figure \ref{fig:download_mfn}.}
\label{alg:download_mfn}
\begin{verbatim}
import numpy as np
import matplotlib.pyplot as plt
plt.rcParams['pgf.texsystem'] = 'pdflatex'

def d_MFN_day(L, n):
    return (n*L**2/(2**L-1) + L)*86400*250/1024**2

n_1 = 1700
n_2 = 10000
n_3 = 50000
L = 30
x = range(1, L)

y1 = [np.log10(d_MFN_day(i, n_1)) for i in x]
y2 = [np.log10(d_MFN_day(i, n_2)) for i in x]
y3 = [np.log10(d_MFN_day(i, n_3)) for i in x]

fig, ax = plt.subplots()
ax.plot(x, y1, label = f'$\lg(d_{{MFN}}(day)), n = {n_1}$',\
    linestyle = '-')
ax.plot(x, y2, label = f'$\lg(d_{{MFN}}(day)), n = {n_2}$',\
    linestyle = '--')
ax.plot(x, y3, label = f'$\lg(d_{{MFN}}(day)), n = {n_3}$',\
    linestyle = '-.')
ax.legend()
plt.xlabel('L')
plt.ylabel('$\lg(\mathrm{MB})$')

plt.savefig('download_mfn.pgf')
	\end{verbatim}
\end{algorithm}

\begin{algorithm}
\caption{Algorithm to plot graph in Figure \ref{fig:optimal_L_mfn}.}
\label{alg:optimal_L_mfn}
\begin{verbatim}
import numpy as np
from scipy import optimize
import matplotlib.pyplot as plt
plt.rcParams['pgf.texsystem']= 'pdflatex'

def s_MFN_day(L, n):
    return n*L/(2**L-1)*86400*250/1024**2

def d_MFN(L, n):
    return n*L**2/(2**L-1) + L

n_range = range(100, 50000, 100)
res = [optimize.minimize(d_MFN, 10, args = [n])\
    for n in n_range]
suc = [(r['success']) for r in res]
if all(suc):
    print('OK!')

L_argmin = [L['x'][0] for L in res]
s_mfn = [s_MFN_day(L, n) for L, n in zip(L_argmin, n_range)]

fig, ax = plt.subplots(2, 1)    
ax[0].plot(n_range, L_argmin, label =\
    '$\mathrm{argmin}(d_{MNF}(L, n))$', linestyle = '-')
ax[1].plot(n_range, s_mfn, label = '$s_{MNF}(day)$'\
    , linestyle = '-')
ax[0].legend()
ax[1].legend()
ax[1].set_xlabel('n')
ax[0].set_ylabel('L')
ax[1].set_ylabel('MB')

plt.savefig('optimal_L_mfn.pgf')
\end{verbatim}
\end{algorithm}

\end{document}